\documentclass[12pt,reqno]{amsart}
\usepackage{mathrsfs}
\usepackage{array}
\usepackage{cases}
\usepackage{epic}
\usepackage{amsfonts}
\usepackage{graphicx}
\usepackage{amsmath}
\usepackage{amssymb, upgreek}
\usepackage{bm}
\usepackage{latexsym,todonotes}
\usepackage{pdflscape}
\usepackage[all]{xypic}
\usepackage[all]{xy}
\usepackage{color}
\usepackage{colordvi}
\usepackage{multicol}
\usepackage[normalem]{ulem}
\usepackage[mathscr]{euscript}
\input diagxy
\usepackage[linktocpage=true]{hyperref}
\hypersetup{colorlinks,linkcolor=blue,urlcolor=cyan,citecolor=blue}
\usepackage{chngcntr} 
\counterwithin{figure}{section}
\counterwithin{table}{section}
\usetikzlibrary{
	cd,
	shapes,
	arrows,
	positioning,
	decorations.markings,
	decorations.pathmorphing,
	circuits.logic.US,
	circuits.logic.IEC,
	fit,
	calc,
	plotmarks,
	matrix
}
\begin{document}
	\input xy
	\xyoption{all}

	\newtheorem{innercustomthm}{{\bf Theorem}}
	\newenvironment{customthm}[1]
	{\renewcommand\theinnercustomthm{#1}\innercustomthm}
	{\endinnercustomthm}
	
	\newtheorem{innercustomcor}{{\bf Corollary}}
	\newenvironment{customcor}[1]
	{\renewcommand\theinnercustomcor{#1}\innercustomcor}
	{\endinnercustomthm}
	
	\newtheorem{innercustomprop}{{\bf Proposition}}
	\newenvironment{customprop}[1]
	{\renewcommand\theinnercustomprop{#1}\innercustomprop}
	{\endinnercustomthm}

	\newtheorem{theorem}{Theorem}[section]
	\newtheorem{acknowledgement}[theorem]{Acknowledgement}
	\newtheorem{algorithm}[theorem]{Algorithm}
	\newtheorem{axiom}[theorem]{Axiom}
	\newtheorem{case}[theorem]{Case}
	\newtheorem{claim}[theorem]{Claim}
	\newtheorem{conclusion}[theorem]{Conclusion}
	\newtheorem{condition}[theorem]{Condition}
	\newtheorem{conjecture}[theorem]{Conjecture}
	\newtheorem{construction}[theorem]{Construction}
	\newtheorem{corollary}[theorem]{Corollary}
	\newtheorem{Claim}[theorem]{Claim}
	\newtheorem{definition}[theorem]{Definition}
	\newtheorem{example}[theorem]{Example}
	\newtheorem{exercise}[theorem]{Exercise}
	\newtheorem{lemma}[theorem]{Lemma}
	\newtheorem{notation}[theorem]{Notation}
	\newtheorem{problem}[theorem]{Problem}
	\newtheorem{proposition}[theorem]{Proposition}
	\newtheorem{solution}[theorem]{Solution}
	\newtheorem{summary}[theorem]{Summary}
	\numberwithin{equation}{section}
	
	\theoremstyle{remark}
	\newtheorem{remark}[theorem]{Remark}
	
	\makeatletter
	\newcommand{\rmnum}[1]{\romannumeral #1}
	\newcommand{\Rmnum}[1]{\expandafter\@slowromancap\romannumeral #1@}
	\def \g{\mathfrak{g}}
	\def \bH{{\mathbf H}}
	\def \bB{{\mathbf B }}
	\def \BKH{{\acute{H}}}
	\def \bBKH{{\acute{\mathbf H}}}
	\def \bK{\mathbb{K}}
	\def \bC{\mathbb{\bK_\de}}
	\def \bQ{\mathbb{Q}}
	\def \Z{\mathbb{Z}}
	\def \I{\mathbb{I}}
	\def \ch{{\mathcal H}}
	\def \cm{{\mathcal M}}
	\def \ct{{\mathcal T}}
	\def \hW{\widehat{W}}
	\def \bZ{\mathbb{Z}}
	\def \TH{\Theta}
	\def \tMH{{\cm\widetilde{\ch}(\Lambda^\imath)}}
	\def \tMHl{{\cm\widetilde{\ch}(\bs_\ell\Lambda^\imath)}}
	\def \bTH { \boldsymbol{\Theta}}
	\def \bDel{ \boldsymbol{\Delta}}
	\def \tY{\widetilde{Y}}
	\def \tH{{\widetilde{H}}}
	\def \tTH{\widetilde{\Theta}}
	\def \btau{\widehat{\tau}}
	\def \s{\varsigma}
	\def \bvs{\boldsymbol{\varsigma}}
	\def \bs{\mathbf r}
	\def \cR{\mathcal{R}}
	\def \fg{\mathfrak{g}}
	\def \Aut{\operatorname{Aut}\nolimits}
	
	\newcommand{\hgt}{\text{ht}}
	\def \ad{\text{ad}\,}
	\def \gr{\text{gr}\,}
	\def \hg{\widehat{\g}}
	\def \mf{\mathfrak}
	\def \N{\mathbb N}
	\def \tk{\widetilde{k}}
	\def \brW{{\rm Br}(W_{\tau})}
	\def \bome{{\varpi}}
	\def \bth{\bm{\theta}}
	\def\Br{{\rm Br}}
	\renewcommand{\t}{\boldsymbol{\omega}}
	
	\def \tfX{\widetilde{ \Upsilon}}
	\def \tT{\widetilde{T}}

	\def \tcT{\widetilde{\mathscr T}}
	\def \E{\widetilde{E}}
	
	\def \SS{\mathbb{S}}
	\def \sll{\mathfrak{sl}}
	\def \tUD{{}^{Dr}\tU}

	\newcommand{\UU}{{\mathbf U}\otimes {\mathbf U}}
	\newcommand{\UUi}{(\UU)^\imath}
	\newcommand{\tUU}{{\tU}\otimes {\tU}}
	\newcommand{\tUUi}{(\tUU)^\imath}
	\newcommand{\tK}{\widetilde{K}}
	\newcommand{\tU}{\widetilde{{\mathbf U}} }
	\newcommand{\tUi}{\widetilde{{\mathbf U}}^\imath}
	\newcommand{\tUiii}{\tUi(\widehat\sll_3,\tau)}
	\newcommand{\tUiD}{{}^{\text{Dr}}\widetilde{{\mathbf U}}^\imath}
	\newcommand{\tUigr}{{\text{gr}}\tUi}
	\newcommand{\tUigrp}{{\text{gr}}\widetilde{{\mathbf U}}^{\imath,+}}
	\newcommand{\tUiDgr}{{\text{gr}}\tUiD}

	\newcommand{\sqq}{{\bf v}}
	\newcommand{\sqvs}{\sqrt{\vs}}
	\newcommand{\dbl}{\operatorname{dbl}\nolimits}
	\newcommand{\swa}{\operatorname{swap}\nolimits}
	\newcommand{\Gp}{\operatorname{Gp}\nolimits}
	\newcommand{\Sym}{\operatorname{Sym}\nolimits}
	\newcommand{\qbinom}[2]{\begin{bmatrix} #1\\#2 \end{bmatrix} }
	
	\def \ov{\overline}
	\def \balpha{{\boldsymbol \alpha}} 
	\def \K{\mathbb K}
	
	\newcommand{\nc}{\newcommand}
	\nc{\browntext}[1]{\textcolor{brown}{#1}}
	\nc{\greentext}[1]{\textcolor{green}{#1}}
	\nc{\redtext}[1]{\textcolor{red}{#1}}
	\nc{\bluetext}[1]{\textcolor{blue}{#1}}
	\nc{\brown}[1]{\browntext{ #1}}
	\nc{\green}[1]{\greentext{ #1}}
	\nc{\red}[1]{\redtext{#1}}
	\nc{\blue}[1]{\bluetext{#1}}

	\def \Q {\mathbb Q}
	\def \C{\mathbb C}
	\def \TT{\mathbf T}
	\def \tt{\mathbf t}
	\newcommand{\wt}{\text{wt}}
	\def \de{\delta}
	\def \bvs{{\boldsymbol{\varsigma}}}
	\def \vs{\varsigma}
	\def \U{\mathbf U}
	\def \Ui{\mathbf{U}^\imath}
	\def \dm{\diamond}
	\def \bS{\mathbb S}
	
	\def \bR{\mathbb R}
	\def \bP{\mathbb P}
	\def \bF{\mathbb F}
	\def \II{\mathbb{I}_0}
    \def \t{\Omega}
	\allowdisplaybreaks
	
	\newcommand{\lutodo}{\todo[inline,color=violet!20, caption={}]}
	\newcommand{\wtodo}{\todo[inline,color=cyan!20, caption={}]}
	\newcommand{\ztodo}{\todo[inline,color=green!20, caption={}]}
	
	\title[Quasi-split affine $\imath$quantum groups II]{Braid group action and quasi-split affine $\imath$quantum groups II: higher rank}
	
	\author[Ming Lu]{Ming Lu}
	\address{Department of Mathematics, Sichuan University, Chengdu 610064, P.R.China}
	\email{luming@scu.edu.cn}

	\author[Weiqiang Wang]{Weiqiang Wang}
	\address{Department of Mathematics, University of Virginia, Charlottesville, VA 22904, USA}
	\email{ww9c@virginia.edu}

	\author[Weinan Zhang]{Weinan Zhang}
	\address{Department of Mathematics and New Cornerstone
Science Laboratory, The University of Hong Kong, Hong Kong SAR, P.R.China}
	\email{mathzwn@hku.hk}


	\begin{abstract}
		This paper studies quantum symmetric pairs $(\widetilde{\mathbf U}, \widetilde{{\mathbf U}}^\imath )$ associated with quasi-split Satake diagrams of affine type $A_{2r-1}, D_r, E_{6}$ with a nontrivial diagram involution fixing the affine simple node. Various real and imaginary root vectors for the universal $\imath$quantum groups $\widetilde{{\mathbf U}}^\imath$ are constructed with the help of the relative braid group action, and they are used to construct affine rank one subalgebras of  $\widetilde{{\mathbf U}}^\imath$. We then establish relations among real and imaginary root vectors in different affine rank one subalgebras and use them to give a Drinfeld type presentation of $\widetilde{{\mathbf U}}^\imath$.
	\end{abstract}
	
	\maketitle

 \begin{quote}
{\em Dedicated to Vyjayanthi Chari for her 65th birthday with admiration
}
\end{quote}

 \setcounter{tocdepth}{1}
	\tableofcontents

	\section{Introduction}

\subsection{}
Associated with Satake diagrams $(\I =\I_{\bullet} \cup \I_{\circ}, \tau)$, G. Letzter \cite{Let02} constructed quantum symmetric pairs $(\U, \Ui)$, where the $\imath$quantum group $\Ui$ is a coideal subalgebra of a Drinfeld-Jimbo quantum group $\U$ (also cf. \cite{Ko14}). The notion of universal $\imath$quantum groups $\tUi$ was introduced by two of the authors \cite{LW22a} owing to their realization through the $\imath$Hall algebra, where $\tUi$ is a coideal subalgebra of a Drinfeld double quantum group $\tU$; Letzter-Kolb's version of $\imath$quantum groups with parameters $\Ui$ can then be recovered by a central reduction from $\tUi$. We work mostly with $\tUi$ in this paper as it provides the most natural setting for the relative braid group symmetries; see \cite{LW22b, WZ23, Z23} and compare \cite{KP11}. 

An affine quantum group $\U$ possesses a Drinfeld (loop) presentation that manifestly exhibits the quantized loop algebra structure of $\U$ \cite{Dr88}. In the Drinfeld presentation for $\U$ of untwisted type, the affine rank one subalgebra $\U(\widehat\sll_2)$ is ubiquitous. Lusztig's braid group action on $\U$ \cite{Lus93} has also played a fundamental role in the construction of real and imaginary root vectors and then in the proof of an isomorphism between the Drinfeld presentation and the Drinfeld-Jimbo presentation of $\U$ \cite{Be94, Da12, Da15}. The Drinfeld presentation was instrumental for the finite-dimensional representation theory of $\U$; cf. \cite{CP91}.

The quasi-split $\imath$quantum groups, which correspond to Satake diagrams with $\I_{\bullet} =\emptyset$, already form a rich family of algebras which generalize Drinfeld-Jimbo quantum groups. In contrast to the untwised affine quantum group setting, there are three quasi-split affine rank one $\imath$quantum groups (assuming that $\tau$ fixes the affine simple node). For one, the quantum group $\U(\widehat\sll_2)$ can be viewed as a rank one quasi-split affine $\imath$quantum group. The second one, the {\em split} affine rank one $\imath$quantum group, is known as the q-Onsager algebra; cf., e.g., \cite[Example 7.6]{Ko14}. The Drinfeld type presentations of the (universal) $\imath$quantum groups of {\em split} ADE type were initiated by two of the authors in \cite{LW21b} and subsequently extended to split BCFG types with streamlined proofs by the third author \cite{Z22}; here ``split" means that $\tau=\text{Id}$ in quasi-split Satake diagrams. A Drinfeld type presentation in the third rank one quasi-split affine $\imath$quantum group was obtained in \cite{LWZ23}. All these constructions are set up in the framework of universal $\imath$quantum groups $\tUi$ where the relative braid group symmetries of $\tUi$ have played a crucial role. 

The study of affine $\imath$quantum groups and twisted Yangians is very relevant to quantum integrable models with boundary conditions (cf. \cite{EKP07}) such as open spin chains of XXZ type \cite{Skl88} and boundary affine Toda field theories (cf. \cite{BB10}). One may hope the Drinfeld presentations developed here may have physical applications in some of these directions.

\subsection{}
The goal of this paper is to provide a Drinfeld type presentation for the universal $\imath$quantum groups associated with quasi-split Satake diagrams of affine type $A_{2r-1}, D_r, E_{6}$ with a nontrivial diagram involution $\tau$ fixing the affine simple node; these Satake diagrams will be denoted by AIII$_{2r-1}^{(\tau)}$, DI$_r^{(\tau)}$, and EI$_{6}^{(\tau)}$ (see Table \ref{tab:Satakediag}), analogous to the notation for Satake diagrams of finite type.

We start by setting up some basics on the relative affine Weyl groups for the underlying affine symmetric pairs $(\widehat\fg, \widehat\fg^{\omega\tau})$, where $\omega$ denotes the Chevalley involution.
Building on our previous constructions \cite{WZ23} (also see \cite{LW21a, Z23}), we construct the relative braid group symmetries of $\tUi$ and develop their basic properties which we need.

Using the relative braid group symmetries, we formulate affine rank one subalgebras $\tUi_{[i]}$ of  $\widetilde{{\mathbf U}}^\imath$ associated with each $\tau$-orbit $\{i,\tau i\}$ of $\I_0$. We construct an algebra homomorphism $\aleph_i$ from a suitable affine rank one $\imath$quantum group into $\tUi$ for each $i\in \I_0$.  In the setting of this paper, there are 2 types of affine rank one subalgebras depending on whether $c_{i,\tau i}=2$ or $0$; that is, the q-Onsager algebra when $c_{i,\tau i}=2$ and (a universal variant of) the affine quantum group $\U(\widehat\sll_2)$ when $c_{i,\tau i}=0$.

It takes some serious effort to show that $\aleph_i$ is an isomorphism onto the rank one subalgebra $\tUi_{[i]}$ of $\tUi$ (see Propositions~\ref{prop:rank1iso} and \ref{prop:rank1isoQG}). While the formulation of these embeddings here is similar to those constructed in \cite{Be94} and \cite{LW21b}, the old arguments do not work for the new generalities of  this paper. In the proof that $\aleph_i$ is a homomorphism, we apply some of new technical tools developed more recently in \cite{WZ23} together with some weight arguments to establish several crucial Serre-type relations involving relative braid group symmetries.

We show that the relative braid group actions are compatible with such embeddings of affine rank one algebras. This is a direct consequence of the rank one formulas for braid group actions on modules \cite{WZ24}. (These rank one formulas imply the compatibility of relative braid group actions under rank one embeddings in all quasi-split types including the split types in \cite{LW21b, Z22}.) This compatibility then allows us to construct real and imaginary root vectors for $\tUi$ by transporting their counterparts in the rank one settings through the rank one embeddings. These root vector constructions are natural generalizations of root vectors for affine quantum groups \cite{Da93, Be94} and for (quasi-)split affine $\imath$quantum groups  \cite{BK20, LW21b, Z22, LWZ23}. 

We establish relations among several families of real and imaginary root vectors and use them to give a Drinfeld type presentation of $\widetilde{{\mathbf U}}^\imath$ (see Definition~\ref{def:iDRA1} and Theorem~\ref{thm:Dr}). Verifying these relations is the most technical part of the paper and occupies the whole Section~ \ref{sec:relation}. Some of the relations can be verified as in the split type, and in these cases, we refer to the arguments given in \cite{LW21b, Z22} to avoid repetition. There are also several additional relations which require new arguments and we give full details. 

A generating function formulation of the current relations and Drinfeld presentation of $\tUi$ is also provided (see Theorem~\ref{thm:DrGF}). We also formulate an alternative Drinfeld presentation of $\tUi$ in Theorem \ref{thm:Drvariant}, which uses a variant of imaginary root vectors. 

\subsection{}
\label{subsec:future}

The Drinfeld presentation of quasi-split $\imath$quantum groups is expected to have many applications, some of which are indicated below. 

The Drinfeld presentation should be instrumental for developing a theory of q-characters for finite-dimensional representations of affine $\imath$quantum groups as the Drinfeld presentation manifestly exhibits a maximally commutative subalgebra of $\tUi$. 

The Drinfeld presentation of affine $\imath$quantum groups will lead via a degeneration approach to a Drinfeld presentation for the twisted Yangians of the corresponding quasi-split types. 

The Drinfeld presentation will be crucially used in the categorical realization of affine $\imath$quantum groups via $\imath$Hall algebras of weighted projective lines coupled with an involution. 

Yet in another direction, the Drinfeld presentation of quasi-split $\imath$quantum groups of affine type AIII will be needed in the geometric realization via equivariant K-groups of Steinberg varieties of classical type. 

\subsection{}
The Drinfeld presentation constructions in this paper (together with \cite{LW21b, Z22}) have covered all $\imath$quantum groups associated with quasi-split affine Satake diagrams with a nontrivial diagram involution $\tau$ fixing the affine simple node, except the type AIII$_{2r}^{(\tau)}$. The special case $r=1$  for AIII$_{2r}^{(\tau)}$ was treated in \cite{LWZ23}, while the general case for $r\geq 2$ remains open. This class for $r\ge 2$ is the only class of quasi-split affine $\imath$quantum groups $\tUi$ with 3 different lengths of relative roots, and its relative braid group is of type $A_{2r}^{(2)}$. We have some partial success on Drinfeld presentation of $\tUi$ of type $A_{2r}^{(\tau)}$ and hope to return to it elsewhere. In this case, the constructions of root vectors are very different and establishing the relations among root vectors is more challenging. (This is reminiscent of the Drinfeld presentation of twisted affine quantum groups of type $A_{2r}^{(2)}$, which needed to be treated separately \cite{Da12, Da15}.) 

The Drinfeld presentation remains to be understood for a quasi-split affine $\imath$quantum group associated to a Satake diagram whose affine simple node is not fixed by the diagram involution $\tau$; the classical limits already look messy in these cases. For the Drinfeld presentations for non-quasi-split affine $\imath$quantum groups, the relative braid group action may be insufficient for constructing the desired real root vectors, though the relative braid group symmetries constructed in \cite{WZ23} will still be useful. 

So far, we have only addressed affine $\imath$quantum groups whose underlying Dynkin diagrams are of untwisted affine types. It will be yet another great challenge to extend the Drinfeld presentations to twisted affine $\imath$quantum groups (see \cite{Da12, Da15} for Drinfeld presentations of twisted affine quantum groups).

\subsection{}
The paper is organized as follows. In Section~\ref{sec:prelim}, we introduce the quasi-split affine $\imath$quantum groups $\tUi$, the related relative affine Weyl groups and affine braid groups. In  Section~\ref{sec:subalgebras}, we construct affine rank one subalgebras of $\tUi$, for each $i\in \II$, and show that they are isomorphic to either $q$-Onsager algebra (if $\tau i=i$) or quantum $\widehat{\sll}_2$ (if $\tau i \neq i$). The proofs of two technical lemmas are postponed to Appendix \ref{app:A1}--\ref{app:A2}. 
In Section~ \ref{sec:Dpresentation}, we use the affine rank one subalgebras to construct various root vectors for $\tUi$, and then formulate the Drinfeld presentation for $\tUi$. The tedious verification of the relations in the Drinfeld presentation is completed in Section \ref{sec:relation}.

	\section{Affine $\imath$quantum groups and relative braid group action}
   \label{sec:prelim}
   
	In this section, we review affine quantum groups, their braid group actions and Drinfeld's presentations. Then we set up  $\imath$quantum groups $\tUi$ of quasi-split affine type $A_{2r-1}, D_r, E_{6}$. We formulate explicitly the corresponding (extended) relative affine Weyl groups, and establish some basic structures of the relative braid group symmetries of $\tUi$.

	\subsection{Affine Weyl and braid groups}

	Let $(c_{ij})_{i,j\in \II}$ be the Cartan matrix of the simple Lie algebra $\fg$ of type ADE. Let $\cR_0$ be the set of roots for $\fg$, and fix a set $\cR^+_0$ of positive roots  with simple roots $\alpha_i$ $(i\in \II)$. Denote by $\theta$ the highest root of $\fg$.
	
	Let $\widehat{\fg}$ be the (untwisted) affine Lie algebra with affine Cartan matrix denoted by $(c_{ij})_{i,j\in\I}$, where $\I=\{0\} \cup \II$ with the affine node $0$. Let $\alpha_i$ $(i\in \I)$ be the simple roots of $\widehat{\fg}$, and $\alpha_0=\de -\theta$, where $\de$ denotes the basic imaginary root. The root system $\cR$  for $\widehat{\fg}$ and its positive system $\cR^+$ are defined to be
	\begin{align}
		\cR &=\{\pm (\beta + k \delta) \mid \beta \in \cR_0^+, k  \in \Z\}  \cup \{m \delta \mid m \in \Z\backslash \{0\} \},
		\label{eq:roots}  \\
		\cR^+ &= \{k \delta +\beta \mid \beta \in \cR_0^+, k  \ge 0\}
		\cup  \{k \delta -\beta \mid \beta \in \cR_0^+, k > 0\}
		\cup \{m \delta \mid m \ge 1\}.
		\label{eq:roots+}
	\end{align}
	For $\gamma =\sum_{i\in \I} n_i \alpha_i \in \N \I$, the height $\text{ht} (\gamma)$ is defined as $\text{ht} (\gamma) =\sum_{i\in \I} n_i$. We refer to \cite{K90} for more details.
	
	Let $P$ and $Q$ denote the weight and root lattices of the simple Lie algebra $\fg$, respectively. Let $\omega_i \in P$ $(i\in \II)$ be the fundamental weights of $\fg$. Note $\alpha_i =\sum_{j\in \II} c_{ij}\omega_j$. We define a bilinear pairing $\langle \cdot, \cdot \rangle : P\times Q \rightarrow \Z$ such that $\langle \omega_i, \alpha_j \rangle =\delta_{ij}$, for $i,j \in \II$, and thus $\langle \alpha_i, \alpha_j \rangle = c_{ij}$.
	
	The Weyl group $W_0$ of $\fg$ is generated by the simple reflections $s_i$, for $i \in \II$. They act on $P$ by $s_i(x)=x-\langle x, \alpha_i \rangle\alpha_i$ for $x\in P$. The extended affine Weyl group $\widetilde{W}$ is the semi-direct product $W_0 \ltimes P$, which contains the affine Weyl group $W:=W_0 \ltimes Q =\langle s_i \mid i \in \I \rangle$ as a subgroup; we denote
	\[
	t_\omega =(1, \omega) \in \widetilde W, \quad \text{ for } \omega \in P.
	\]
	In particular, for $\omega\in P$, $j\in\II$, $t_\omega(\alpha_j)=\alpha_j-\langle \omega,\alpha_j\rangle \de$.
	We identify $P/Q$ with a finite group $\Omega$ of Dynkin diagram automorphism, and so $\widetilde{W} =\Omega . W$. There is a length function $\ell(\cdot)$ on $\widetilde{W}$ such that $\ell(s_i)=1$, for $i\in \I$, and each element in $\Omega$ has length 0.
	
	For $i\in \II$, as in \cite{Be94}, we have
	\begin{equation}  \label{eq:tomega}
		\ell(t_{\omega_i}) =\ell(t_{\omega_i'})+1, \qquad
		\text{ where }  t_{\omega_i'}:= t_{\omega_i} s_i \in W.
	\end{equation}
	
	Let $\tU =\tU(\widehat{\fg})$ denote the Drinfeld double affine quantum group, a $\Q(v)$-algebra generated by $E_i, F_i, \tK_i,\tK_i'$, for $i\in \I$ subject to the following relations:
\begin{align}
[E_i,F_j]= \delta_{ij} \frac{\tK_i-\tK_i'}{v-v^{-1}},  &\qquad [\tK_i,\tK_j]=[\tK_i,\tK_j']  =[\tK_i',\tK_j']=0,
\label{eq:KK}
\\
\tK_i E_j=v^{c_{ij}} E_j \tK_i, & \qquad \tK_i F_j=v^{-c_{ij}} F_j \tK_i,
\label{eq:EK}
\\
\tK_i' E_j=v^{-c_{ij}} E_j \tK_i', & \qquad \tK_i' F_j=v^{c_{ij}} F_j \tK_i',
 \label{eq:K2}
\end{align}
and the quantum Serre relations (which we skip). Here $\tK_i\tK_i'$ are central in $\tU$, for all $i\in\I$. 
Then the Drinfeld-Jimbo quantum group $\U$ \cite{Dr87,J85} for $\widehat{\fg}$ is obtained by the central reduction from $\tU$:
$$\U=\tU/(\tK_i\tK_i'-1\mid i\in\I). $$
Denote the $v$-integers by $[n] =\frac{v^n-v^{-n}}{v-v^{-1}}$ and recall the divided powers 
\[
F_i^{(n)} =F_i^n/[n]^!, \quad E_i^{(n)} =E_i^n/[n]^!, \quad \text{ for } n\ge 1 \text{ and  } i\in \I.
\]
Following Lusztig \cite{Lus93}, we have the braid group action on $\U$, which can be easily adapted to a braid group action on $\tU$; see \cite[Proposition 6.20]{LW22b}.
	
	\begin{lemma}
		\label{prop:BG1U}
		For $i\in \I$, there are automorphisms $T_i$ on $\tU$  such that
		\begin{align*}
			T_i(\tK_j)&= \tK_i^{-c_{ij}} \tK_j,\qquad T_i(\tK'_j)= (\tK'_i)^{-c_{ij}} \tK'_j,
			\\
			T_{i}(F_i) &=  -\tK_{i}^{-1}E_{ i}, \qquad T_{i}(E_i)= -F_i(\tK'_{i})^{-1},
			\\
			T_{i}(F_j) &=  \sum^{-c_{ i,j}}_{r=0} (-1)^{r} v^{ r} F_i^{(r)}  F_j F_i^{(-c_{ij}-r)}\qquad (j\neq i),
			\\
			T_{i}(E_j)&=  \sum^{-c_{ i,j}}_{r=0}  (-1)^{r} v^{-r} E_{i}^{(-c_{\tau i,j}-r)} E_j E_{i}^{(r)}\qquad (j\neq i).
		\end{align*}
	\end{lemma}
 These automorphisms $T_i$ descend to automorphisms of $\U$, also denoted by $T_i$. It is well known that they give rise to automorphisms $T_w$ on $\tU$ and $\U$, for $w \in \widetilde W$.
	
	The following was known to Bernstein and Lusztig.
	\begin{lemma} [\text{\cite[\S2.7]{Lus89}}]
		\label{lem:Braid-Lus}
		Let $x\in P$, $i, j \in \II$.
		\begin{enumerate}
			\item[(a)]
			If $s_ix=xs_i$, then $T_iT_x=T_xT_i$.
			\item[(b)]
			If $s_i xs_i=t_{\alpha_i}^{-1}x=\prod_{k\in \II} t_{\omega_k}^{a_k}$, then we have $T_i^{-1}T_xT_i^{-1}=\prod_k T_{\omega_k}^{a_k}$; in particular, we have $T_i^{-1}T_{\omega_i}T_i^{-1}=T_{\omega_i}^{-1}\prod_{k\neq i}T_{\omega_k}^{-c_{ik}}$.
			\item[(c)]
			$T_{\omega_i} T_{\omega_j}  = T_{\omega_j} T_{\omega_i}.$
		\end{enumerate}
	\end{lemma}

	\subsection{Drinfeld presentation of affine quantum groups}
	\label{subsec:Drsl2}
 
	The affine quantum group $\U$ admits a second presentation known as the Drinfeld presentation \cite{Dr88}. Recall $(c_{ij})_{i,j\in\I_0}$ is the Cartan matrix of the simple Lie algebra $\g$. Let $^{\text{Dr}}\U$ be the $\Q(v)$-algebra generated by $x_{i k}^{\pm}$, $h_{i l}$, $K_i^{\pm 1}$, $\texttt{C}^{\pm \frac12}$, for $i\in\I_0$, $k\in\Z$, and $l\in\Z\backslash\{0\}$, subject to the following relations: $\texttt{C}^{\frac12}, \texttt{C}^{- \frac12}$ are central and
	\begin{align*}
		[K_i,K_j] & =  [K_i,h_{j l}] =0, \quad K_i K_i^{-1} =\texttt{C}^{\frac12} \texttt{C}^{- \frac12} =1,
		\\
		[h_{ik},h_{jl}] &= \delta_{k, -l} \frac{[k c_{ij}]}{k} \frac{\texttt{C}^k -\texttt{C}^{-k}}{v -v^{-1}},
		\\
		K_ix_{jk}^{\pm} K_i^{-1} &=v^{\pm c_{ij}} x_{jk}^{\pm},
		\\
		[h_{i k},x_{j l}^{\pm}] &=\pm\frac{[kc_{ij}]}{k} \texttt{C}^{\mp \frac{|k|}2} x_{j,k+l}^{\pm},
		\\
		[x_{i k}^+,x_{j l}^-] &=\delta_{ij} {(\texttt{C}^{\frac{k-l}2} K_i\psi_{i,k+l} - \texttt{C}^{\frac{l-k}2} K_i^{-1} \varphi_{i,k+l})}, 
		\\
		x_{i,k+1}^{\pm} x_{j,l}^{\pm}-v^{\pm c_{ij}} x_{j,l}^{\pm} x_{i,k+1}^{\pm} &=v^{\pm c_{ij}} x_{i,k}^{\pm} x_{j,l+1}^{\pm}- x_{j,l+1}^{\pm} x_{i,k}^{\pm},
		\\
		\Sym_{k_1,\dots,k_r}\sum_{t=0}^{r} (-1)^t \qbinom{r}{t}
		& x_{i,k_1}^{\pm}\cdots
		x_{i,k_t}^{\pm} x_{j,l}^{\pm}  x_{i,k_t+1}^{\pm} \cdots x_{i,k_n}^{\pm} =0, \text{ for } r= 1-c_{ij}\; (i\neq j),
	\end{align*}
	where
	$\Sym_{k_1,\dots,k_r}$ denotes the symmetrization with respect to the indices $k_1,\dots,k_r$, $\psi_{i,k}$ and $\varphi_{i,k}$ are defined by the following equations:
	\begin{align*}
		1+ \sum_{m\geq 1} (v-v^{-1})\psi_{i,m}u^m &=  \exp\Big((v -v^{-1}) \sum_{m\ge 1}  h_{i,m}u^m\Big),
		\\
		1+ \sum_{m\geq1 } (v-v^{-1}) \varphi_{i, -m}u^{-m} &= \exp \Big((v -v^{-1}) \sum_{m\ge 1} h_{i,-m}u^{-m}\Big).
	\end{align*}
	(We omit a degree operator $D$ in the version of ${}^{\text{Dr}}\U$ above.)
	
	There exists an isomorphism of $\Q(v)$-algebras \cite{Dr88, Be94, Da15}
	\begin{align}
		\label{eq:phi1}
		\phi:  {}^{\text{Dr}}\U \longrightarrow \U.
	\end{align}

	\subsection{(Universal) $\imath$quantum groups}
	
	For a (generalized) Cartan matrix $C=(c_{ij})_{\I\times\I}$, let $\Aut(C)$ be the group of all permutations $\tau$ of the set $\I$ such that $c_{ij}=c_{\tau i,\tau j}$. 	Let $\tau$ be an involution in $\Aut(C)$, i.e., $\tau^2={\rm Id}$. Following \cite{LW22a} we define $\widetilde{\U}^\imath$ to be the $\Q(v)$-subalgebra of $\tU$ generated by
	\begin{equation}
		\label{eq:Bi}
		B_i= F_i +  E_{\tau i} \tK_i',
		\qquad \tk_i = \tK_i \tK_{\tau i}', \quad \forall i \in \I.
	\end{equation}
	Let $\tU^{\imath 0}$ be the $\Q(v)$-subalgebra of $\tUi$ generated by $\tk_i$, for $i\in \I$. By \cite[Lemma 6.1]{LW22a}, the elements $\tk_i$ (for $i= \tau i$) and $\tk_i \tk_{\tau i}$  (for $i\neq \tau i$) are central in $\tUi$.
	
	Let $\bvs=(\vs_i)\in  (\Q(v)^\times)^{\I}$ be such that $\vs_i=\vs_{\tau i}$ for all $i$. We fix
\begin{align}   \label{eq:ci}
\I_\tau &= \{ \text{the chosen representatives of $\tau$-orbits in $\I$} \}.
\end{align}
Let $\Ui:=\Ui_{\bvs}$ be the $\Q(v)$-subalgebra of $\U$ generated by
	\[
	B_i= F_i+\vs_i E_{\tau i}K_i^{-1},
	\quad
	k_j= K_jK_{\tau j}^{-1},
	\qquad  \forall i \in \I, j \in \I\backslash\I_\tau.
	\]
	It is known \cite{Let02, Ko14} that $\U^\imath$ is a right coideal subalgebra of $\U$ in the sense that $\Delta: \Ui \rightarrow \Ui\otimes \U$; and $(\U,\Ui)$ is called a \emph{quantum symmetric pair} ({\em QSP} for short), and they specialize at $v=1$ to $(U(\fg), U(\fg^{\omega\tau}))$.
	The algebras $\Ui_{\bvs}$, for $\bvs \in  (\Q(v)^\times)^{\I}$, are obtained from $\tUi$ by central reductions. Indeed, by \cite[Proposition 6.2]{LW22a}, the algebra $\Ui$ is isomorphic to the quotient of $\tUi$ by the ideal generated by
	\begin{align}   \label{eq:parameters}
		\tk_i - \vs_i \; (\text{for } i =\tau i),
			\qquad  \tk_i \tk_{\tau i} - \vs_i \vs_{\tau i}  \;(\text{for } i \neq \tau i).
	\end{align}
	The isomorphism is given by sending $B_i \mapsto B_i, k_j \mapsto \vs_{\tau j}^{-1} \tk_j, k_j^{-1} \mapsto \vs_{j}^{-1} \tk_{\tau j}, \forall i\in \I, j\in \I\backslash\I_\tau$. 

	We shall refer to $\tUi$ and $\Ui$ as {\em (quasi-split) $\imath${}quantum groups}; they are called {\em split} if $\tau =\mathrm{Id}$.
	
	For any $i\in\I$, we set
	\begin{align*}
		\K_i:=-v^{2}\tk_i, \text{ if }\tau i=i;
		\qquad
		\K_j=\tk_j, \text{ otherwise.}
	\end{align*}

	\subsection{Quasi-split affine $\imath$quantum groups}
 \label{subsec:qsAffine}
	
	From now on, let $(c_{ij})_{i,j\in \I}$ be the Cartan matrix of affine ADE type, for $\I = \II \cup \{0\}$ with the affine node $0$. 
	We always assume the involution $\tau$ fixes the affine node $0$ in this paper. 
	
	The universal quasi-split affine $\imath$quantum group admits the following presentation \cite{CLW21}: the $\Q(v)$-algebra $\tUi =\tUi(\widehat{\fg})$ is generated by $B_i$, $\K_i^{\pm 1}$ $(i\in \I)$, subject to the following relations \eqref{eq:SK}--\eqref{eq:S3}, for $i, j\in \I$:
	\begin{align}
		\K_i\K_i^{-1} =\K_i^{-1}\K_i=1,  \quad \K_i\K_\ell &=\K_\ell \K_i, \quad
		\K_\ell B_i=v^{c_{\tau \ell,i} -c_{\ell i}} B_i \K_\ell,
		\label{eq:SK} \\
		B_iB_j -B_j B_i&=0, \qquad\qquad\qquad \text{ if } c_{ij}=0 \text{ and }\tau i\neq j,
		\label{eq:S1} \\
		\sum_{s=0}^{1-c_{ij}} (-1)^s \qbinom{1-c_{ij}}{s} & B_i^{s}B_jB_i^{1-c_{ij}-s} =0, \quad \text{ if } j \neq \tau i\neq i,
		\label{eq:S6}
		\\
		B_{\tau i}B_i -B_i B_{\tau i}& =   \frac{\K_i -\K_{\tau i}}{v-v^{-1}},
		\quad \text{ if }  c_{i,\tau i}=0,
		\label{relation5}	\\
		B_i^2 B_j -[2] B_i B_j B_i +B_j B_i^2 &= - v^{-1}  B_j \K_i,  \qquad \text{ if }c_{ij}=-1 \text{ and }c_{i,\tau i}=2,
		\label{eq:S2}
  \\
		\sum_{r=0}^3 (-1)^r \qbinom{3}{r} B_i^{3-r} B_j B_i^{r} &= -v^{-1} [2]^2 (B_iB_j-B_jB_i) \K_i, 
		\quad \text{ if } -c_{ij}= c_{i,\tau i}=2.  \label{eq:S3}
	\end{align}

\begin{remark}
\label{rem:Onsager}
    The Serre relation \eqref{eq:S3} arises only in the $\imath$quantum group of split affine rank one (also known as the $q$-Onsager algebra). The universal $q$-Onsager algebra $\tUi(\widehat{\mathfrak{sl}}_2)$ is the $\Q(v)$-algebra generated by $B_0,B_1$, $\K_0,\K_1$, subject to the following relations: $\K_0,\K_1$ are central, and
\begin{align}
\sum_{r=0}^3 (-1)^r \qbinom{3}{r} B_i^{3-r} B_j B_i^{r}&= -v^{-1} [2]^2 (B_iB_j-B_jB_i) \K_i, 
\quad \text{ for } i\neq j \in \{0,1\}.
\label{relation:s3}
\end{align}
Below we identify $C =\K_0\K_1$ in $\tUi(\widehat{\mathfrak{sl}}_2)$. We shall need the q-Onsager algebra to construct affine rank one subalgebras of $\tUi$ of higher rank. 
\end{remark}
The split Satake diagrams (i.e., those with $\tau ={\rm Id}$) are identical to the usual Dynkin diagrams. The Drinfeld presentations for the corresponding split affine $\imath$quantum groups were given in \cite{LW21b, Z22}. In this paper, we focus on quasi-split affine $\imath$quantum groups with $\tau\neq {\rm Id}$, all $c_{i,\tau i}\in\{0,2\}$, and the affine node $0$ fixed; accordingly, a Serre relation between $B_i, B_j$ with $c_{i,\tau i}=-1$ does not feature above. A complete list of such quasi-split affine Satake diagrams is given in Table \ref{tab:Satakediag} below, where we have used a superscript $(\tau)$ to denote this class of quasi-split affine Satake diagrams.

 \begin{table}[h]  
\begin{center}
\centering
\begin{tabular}{|m{5cm}<{\centering}|m{8cm}<{\centering}|}
\hline
Types &  Quasi-split affine Satake diagrams \\
\hline

${\rm AIII}_{2r-1}^{(\tau)}$ $(r\geq 2)$ & \setlength{\unitlength}{0.6mm}			\begin{picture}(70,35)(0,5)

   
				\put(0,10){$\circ$}
				\put(0,30){$\circ$}
				\put(50,10){$\circ$}
				\put(50,30){$\circ$}
				\put(72,10){$\circ$}
				\put(72,30){$\circ$}
				\put(92,20){$\circ$}
				\put(-5,5.5){$2r-1$}
				\put(-2,34){${1}$}
				\put(47,6){\small $r+2$}
				\put(47,34){\small $r-2$}
				\put(69,6){\small $r+1$}
				\put(69,34){\small $r-1$}
				\put(92,16){\small $r$}
				
				\put(3,11.5){\line(1,0){16}}
				\put(3,31.5){\line(1,0){16}}
				\put(23,10){$\cdots$}
				\put(23,30){$\cdots$}
				\put(33.5,11.5){\line(1,0){16}}
				\put(33.5,31.5){\line(1,0){16}}
				\put(53,11.5){\line(1,0){18.5}}
				\put(53,31.5){\line(1,0){18.5}}
				\put(-17.5,22.5){\line(2,1){17}}
				\put(-17.5,20.5){\line(2,-1){17}}
				\put(75,12){\line(2,1){17}}
				\put(75,31){\line(2,-1){17}}
				\put(-21,20){$\circ$}
				\put(-21,15){\small $0$}
				
				\color{red}
				\qbezier(0,13.5)(-4,21.5)(0,29.5)
				\put(-0.25,14){\vector(1,-2){0.5}}
				\put(-0.25,29){\vector(1,2){0.5}}
				
				\qbezier(50,13.5)(46,21.5)(50,29.5)
				\put(49.75,14){\vector(1,-2){0.5}}
				\put(49.75,29){\vector(1,2){0.5}}

				\qbezier(72,13.5)(68,21.5)(72,29.5)
				\put(71.75,14){\vector(1,-2){0.5}}
				\put(71.75,29){\vector(1,2){0.5}}
			\end{picture}\\
\hline
${\rm DI}_r^{(\tau)}$ $(r\geq4)$ & \setlength{\unitlength}{0.6mm}	\begin{picture}(40,35)(20,-15)
				\put(0,9.5){$\circ$}
				\put(0,13.5){\small$1$}
				
				\put(3,9.5){\line(2,-1){16.5}}
				\put(3,-8.5){\line(2,1){16.5}}
				\put(0,-11){$\circ$}
				\put(0,-16){\small$0$}
				\put(20,-1){$\circ$}
				\put(20,-6){\small$2$}
				\put(64,-1){$\circ$}
				\put(56,-6){\small$r-2$}
				\put(84,-10){$\circ$}
				\put(80,-13){\small${r-1}$}
				\put(84,9.5){$\circ$}
				\put(84,13.5){\small${r}$}

				\put(38,0){\line(-1,0){15.5}}
				\put(64,0){\line(-1,0){15}}
				
				\put(40,-1){$\cdots$}
				
				\put(83.5,9.5){\line(-2,-1){16.5}}
				\put(83.5,-8.5){\line(-2,1){16.5}}
				
				\put(12,-20.5){\begin{picture}(100,40)	\color{red}
						\qbezier(75,13.5)(79,21.5)(75,29.5)
						\put(75.25,14){\vector(-1,-2){0.5}}
						\put(75.25,29){\vector(-1,2){0.5}}
				\end{picture}}
			\end{picture}
\\
\hline
${\rm EI}_6^{(\tau)}$  &	\setlength{\unitlength}{0.6mm}		\begin{picture}(70,35)(20,7)
				\put(97,36){\small${6}$}
				\put(75,36){\small${5}$}
				\put(97,6.5){\small${1}$}
				\put(75,6.5){\small${2}$}
				\put(11,16){\small $0$}
				\put(33,16){\small $4$}	
				\put(55,16){\small $3$}	\put(10,35){\rotatebox[origin=c]{180}{\begin{picture}(100,10)
							
							\put(10,10){$\circ$}
							
							\put(32,10){$\circ$}
							
							\put(10,30){$\circ$}
							
							\put(32,30){$\circ$}
							\put(96,20){$\circ$}

							\put(76.7,21.2){\line(1,0){19}}
							
							\put(31.5,11){\line(-1,0){19}}
							\put(31.5,31){\line(-1,0){19}}
							
							\put(52,22){\line(-2,1){17.5}}
							\put(52,20){\line(-2,-1){17.5}}
							
							\put(54.7,21.2){\line(1,0){19}}

							\put(52,20){$\circ$}
							
							\put(74,20){$\circ$}
					\end{picture}}

					\put(-37,-32.5){\begin{picture}(100,40)	\color{red}
							\qbezier(5,13.5)(9,21.5)(5,29.5)
							\put(5.25,14){\vector(-1,-2){0.5}}
							\put(5.25,29){\vector(-1,2){0.5}}
					\end{picture}}
					\put(-15,-32.5){\begin{picture}(100,40)	\color{red}
							\qbezier(5,13.5)(9,21.5)(5,29.5)
							\put(5.25,14){\vector(-1,-2){0.5}}
							\put(5.25,29){\vector(-1,2){0.5}}
					\end{picture}}
				}
				
			\end{picture}
\\\hline
\end{tabular} 
\end{center}
\vspace{0.5cm}
\caption{Quasi-split affine Satake diagrams with $\tau \neq {\rm Id}, \tau(0)=0,c_{i,\tau i}\in\{0,2\}$}
\label{tab:Satakediag}
\end{table}
		
For $\mu = \mu' +\mu''  \in \Z \I := \oplus_{i\in \I} \Z \alpha_i$,  define $\K_\mu$ such that
	\begin{align}
		\K_{\alpha_i} =\K_i, \quad
		\K_{-\alpha_i} =\K_i^{-1}, \quad
		\K_{\mu} =\K_{\mu'} \K_{\mu''},
		\quad  \K_\delta =\K_0 \K_\theta.
	\end{align}
	The algebra $\tUi$ is endowed with a filtered algebra structure
	\begin{align}  \label{eq:filt1}
		\widetilde{\U}^{\imath,0} \subset \widetilde{\U}^{\imath,1} \subset \cdots \subset \widetilde{\U}^{\imath,m} \subset \cdots
	\end{align}
	by setting 
	\begin{align}  \label{eq:filt}
		\widetilde{\U}^{\imath,m} =\Q(v)\text{-span} \{ B_{i_1} B_{i_2} \ldots B_{i_s} \K_\mu \mid \mu \in \N\I, i_1, \ldots, i_s \in \I, s \le m \}.
	\end{align}
	Note that
	\begin{align}  \label{eq:UiCartan}
		\widetilde{\U}^{\imath,0} =\bigoplus_{\mu \in \N\I} \Q(v) \K_\mu
	\end{align}
	is the $\Q(v)$-subalgebra generated by $\K_i$ for $i\in \I$. Then, according to a basic result of Letzter and Kolb on quantum symmetric pairs adapted to our setting of $\tUi$ (cf. \cite{Ko14}), the associated graded $\gr \tUi$ with respect to the filtration \eqref{eq:filt1}--\eqref{eq:filt} can be identified as
	\begin{align}   \label{eq:filter}
		\gr \tUi \cong \U^- \otimes \Q(v)[\K_i^\pm | i\in \I],
		\qquad \overline{B_i}\mapsto F_i,  \quad
		\overline{\K}_i \mapsto \K_i \; (i\in \I).
	\end{align}

	\subsection{Relative affine Weyl groups}
	\label{sub:Weyl}

The relative root systems and relative Weyl groups (of finite type) are well known; we refer to the exposition in \cite[\S2.3]{DK19} and the references therein (see \cite{Lus89}). In this subsection we shall adapt this to set up an affine version of relative root systems and relative Weyl groups which are needed in this paper (We do not claim any originality here; see \cite{Lus89}). In contrast, for the split type studied in \cite{LW21b, Z22} we only need the usual affine Weyl groups. 

 Recall $\cR$ is an affine root system. Associated to each $\alpha \in \cR$, we define the element $\balpha\in \Q\cR$ by
	\begin{equation}
		\balpha:=\frac{\alpha+\tau \alpha}{2}, \quad(\alpha\in\cR),
	\end{equation}
and especially we have $\balpha_i$, for $i\in \I$. Note that $\balpha_i=\balpha_{\tau i}$ for $i\in\I$ and $\balpha_0=\alpha_0$. 
 
Let 
\[
\ov{\cR}:=\big\{\balpha =(\alpha+\tau \alpha)/{2}\mid\alpha\in \cR \big\}
\]
be the relative affine root system associated to the quasi-split affine symmetric pair $(\widehat\fg, \widehat\fg^{\omega\tau})$. Let $\I_\tau$ be a fixed set of representatives of $\tau$-orbits on $\I$. Then 
 $\ov{\cR}$ admits a simple system $\{\balpha_i\mid i\in\I_\tau\}$ and the corresponding positive system $\ov{\cR}^+=\{\balpha\mid\alpha\in \cR^+\}$. Let $(\ov{c}_{ij})_{i,j\in \I_\tau}$ be the Cartan matrix of the relative root system, where 
 \[
 \ov{c}_{ij} =\frac{2(\balpha_i,\balpha_j)}{(\balpha_i,\balpha_i)}, \qquad \text{ for } 
 i,j\in \I_\tau.
 \]
	
	We denote by $\bs_{i}$ the following element of order 2 in the Weyl group $W$, i.e.,
	\begin{align}
		\label{def:simple reflection}
		\bs_i= \left\{
		\begin{array}{ll}
			s_{i}, & \text{ if } c_{i,\tau i}= 2 \; (\text{i.e.}, \tau i=i);
			\\
			s_is_{\tau i}, & \text{ if }  c_{i,\tau i}=0;
		\end{array}
		\right.
	\end{align}
	Note that $\bs_i=\bs_{\tau i}$ for any $i\in\I$ and hence we can parametrize the $\bs_i$ by $\I_\tau$.  
 Consider the following subgroup $W_\tau$ of the affine Weyl group $W$ associated with $\cR$:
	\begin{align}
		\label{eq:Wtau}
		W_{\tau} =\{w\in W\mid \tau w =w \tau\}.
	\end{align}
 The group $W_{\tau}$ is a Coxeter group with $\bs_i$ ($i\in \I_\tau$) as its generating set (cf. \cite{Lus03} and \cite[\S2.3]{DK19}), and it is the Weyl group associated to the root system $\ov{\cR}$. We shall refer to $W_{\tau}$ as the {\em relative affine Weyl group} associated with the affine symmetric pair $(\widehat\fg, \widehat\fg^{\omega\tau})$. Note that $W^{\tau}_0 =\{w\in W_0\mid \tau w =w \tau\}$ is the relative (finite) Weyl group associated with the symmetric pair $(\fg, \fg^{\omega\tau})$ with a generating set $\{\bs_i \mid i\in \I_{\tau}\setminus\{0\} \}$. We have the decomposition $W_\tau=W_0^\tau \ltimes Q^\tau$ where $Q^\tau=\{x\in Q\mid\tau (x)=x\}$.

 Similarly, we define the extended relative affine Weyl group $\widetilde{W}_\tau=\{w\in \widetilde{W}\mid \tau w =w \tau\}$ and we have $\widetilde{W}_\tau=W_0^\tau \ltimes P^\tau$ where $P^\tau=\{x\in P\mid\tau (x)=x\}$. We identify $P^\tau /Q^\tau$ with a finite group $\ov{\Gamma}$ of diagram automorphisms on $\ov{\cR}$, and then we have
 \begin{align*}
 \widetilde{W}_\tau\cong W_\tau \rtimes \ov{\Gamma}.
 \end{align*}
 
\begin{remark}
If $\tau$ is trivial, then $W_\tau=W$ and $ \widetilde{W}_\tau= \widetilde{W}$, i.e., (extended) relative affine Weyl groups for split types (considered in \cite{LW21b,Z23}) are just the ordinary (extended) affine Weyl groups.

\end{remark}
 
  Let $\ell^\circ(\cdot)$ be the \emph{length function} of $W_\tau$ (which is extended to a length function on $\widetilde{W}_\tau$ by declaring the length of elements in $\ov{\Gamma}$ to be 0). Define $t_{\bome_i}$ for $i\in \I_\tau\setminus\{0\}$ to be the following elements in the extended Weyl group $\widetilde{W}$:
	\begin{align} \label{tbome2}
		t_{\bome_i}= \left\{
		\begin{array}{ll}
			t_{\omega_i}, & \text{ if } \tau i=i;
			\\
			t_{\omega_i}t_{\omega_{\tau i}}, & \text{ if } \tau i\neq i.
		\end{array}
		\right.
	\end{align}
Then $t_{\bome_i} \in \widetilde{W}_\tau$.

The types of (extended) relative affine Weyl groups associated with the affine Satake diagrams in Table \ref{tab:Satakediag} are specified in Table \ref{tab1}.

\begin{table}[h]  
\begin{center}
\centering
\begin{tabular}{|m{4cm}<{\centering}|m{4cm}<{\centering}|m{6cm}<{\centering}|}
\hline
Types of affine Satake diagrams & Types of relative affine Weyl groups & Dynkin diagrams for relative root systems\\
\hline
${\rm AIII}_{2r-1}^{(\tau)}$ $(r\ge2)$ & ${\rm C}_{r}^{(1)}$ & 		\begin{tikzpicture}[baseline=0, scale=1.5]
			\node at (-0.5,0.4) {$\quad$};
			\node at (-0.5,0.2) {$\circ$};
			\node at (0,0.2) {$\circ$};
			\node at (0.5,0.2) {$\circ$};
			\node at (1.5,0.2) {$\circ$};
			\node at (2.0,0.2) {$\circ$};
			\node at (2.5,0.2) {$\circ$};
			\draw[-] (0.05,0.2) to (0.45,0.2);
			\draw[-] (1.55,0.2) to (1.95,0.2);
			\draw[dashed] (0.55,0.2) to (1.45,0.2);
			\draw[-implies, double equal sign distance]  (2.45, 0.2) to (2.05, 0.2);
			\draw[-implies, double equal sign distance]  (-0.45, 0.2) to (-0.05, 0.2);
			\node at (2,0) {\small $r-1$};
			\node at (2.5,0) {\small $r$};
			\node at (0,0) {\small $1$};
			\node at (0.5,0) {\small $2$};
			\node at (-0.5,0) {\small $0$};
			\node at (-0.5,-0.1) {$\quad$};
		\end{tikzpicture}\\
\hline
${\rm DI}_r^{(\tau)}$ $(r\ge4)$ & ${\rm B}_r^{(1)}$ & \begin{tikzpicture}[baseline=0.0, scale=1.5]
			\node at (-0.5,0.1) {$\circ$};
			\node at (-0.5,0.7) {$\circ$};
			\node at (0,0.4) {$\circ$};
			\node at (0.5,0.4) {$\circ$};
			\node at (1.5,0.4) {$\circ$};
			\node at (2.0,0.4) {$\circ$};
			\node at (2.5,0.4) {$\circ$};
			\draw[-] (0.05,0.4) to (0.45,0.4);
			\draw[-] (-0.45,0.67) to (-0.05,0.44);
			\draw[-] (-0.45,0.14) to (-0.05,0.37);
			\draw[-] (1.55,0.4) to (1.95,0.4);
			\draw[dashed] (0.55,0.4) to (1.45,0.4);
			\draw[-implies, double equal sign distance]  (2.05, 0.4) to (2.45, 0.4);
			\node at (2,0.2) {\small $r-1$};
			\node at (2.5,0.2) {\small $r$};
			\node at (0,0.2) {\small $2$};
			\node at (0.5,0.2) {\small $3$};
			\node at (-0.5,-.1) {\small $1$};
			\node at (-0.5,0.9) {\small $0$};
		\end{tikzpicture}
\\
\hline
${\rm EI}_6^{(\tau)}$ & ${\rm F}_4^{(1)}$ &\begin{tikzpicture}[baseline=0, scale=1.5]
			\node	at (-0.5,0.4){$\quad$};
			\node at (-0.5,0.2) {$\circ$};
			\node at (0,0.2) {$\circ$};
			\node at (0.5,0.2) {$\circ$};
			\node at (1.0,0.2) {$\circ$};
			\node at (1.5,0.2) {$\circ$};
			\draw[-] (-0.45,0.2) to (-0.05,0.2);
			\draw[-] (0.05,0.2) to (0.45,0.2);
			\draw[-] (1.05,0.2) to (1.45,0.2);
		\draw[-implies, double equal sign distance]  (0.55, 0.2) to (0.95, 0.2);
			\node at (1.5,0) {\small $4$};
			\node at (0,0) {\small $1$};
			\node at (0.5,0) {\small $2$};
			\node at (1,0) {\small $3$};
			\node at (-0.5,0) {\small $0$};
			\node	at (-0.5,-0.1){$\quad$};
		\end{tikzpicture}
\\\hline
\end{tabular} 
\end{center}
\vspace{0.5cm}
\caption{Types of relative affine Weyl groups}
\label{tab1}
\end{table}

	\begin{lemma}
		Any reduced expression of $t_{\bome_i}$, for $i\in\II$, ends with $\bs_i$.
	\end{lemma}
	
	\begin{proof}
		If $\ell(t_{\bome_i}\bs_j)<\ell(t_{\bome_i})$, then $t_{\bome_i}(\balpha_j)<0$, which happens if and only if $i=j$.
	\end{proof}

	\subsection{Relative braid group actions}
	
	We shall denote by ${\rm Br}(W)$ the braid group associated to the Weyl group $W$ of $\fg$, and sometimes simply denote by ${\rm Br}(X_k)$ the braid group associated to a Weyl group of type $X_k$. Then the braid group associated to the relative Weyl group for $(\fg, \fg^\theta)$ is of the form
	\begin{equation}
		\label{eq:braidCox}
		\brW =\langle r_i \mid i\in \I_\tau \rangle
	\end{equation}
	where $r_i$ satisfy the same braid relations as for $\bs_i$ in $W_{\tau}$ (but no quadratic relations on $r_i$ are imposed). According to Table \ref{tab1}, the relative braid groups $\brW$ associated with the quasi-split affine Satake diagrams in Table~\ref{tab:Satakediag}  are of type $C_r^{(1)},  B_r^{(1)}, F_4^{(1)}$, respectively. The relative braid group $\brW$ associated with the split affine rank one $\tUi$ (i.e., q-Onsager algebra) is of type $A_1^{(1)}$. The following result was proved in \cite{LW22b} using $\imath$Hall algebra technique (cf. \cite{KP11, BK20, LW21a, WZ23} for earlier formulation in special cases) and a purely algebraic proof was obtained subsequently in \cite{Z23}. 
 
	\begin{theorem} [\text{cf. \cite{LW22b}, \cite[Theorems 3.1, 3.7, 3.11]{Z23}}]
		\label{thm:Ti}
		(1) For $i\in \I$ such that $\tau i=i$, there exists an automorphism $\TT_i$ of the $\Q(v)$-algebra $\tUi$ such that
		$\TT_i(\K_\mu) =\K_{\bs_i\mu}$, and
		\[
		\TT_i(B_j)= \begin{cases}
			\K_i^{-1} B_i,  &\text{ if }j=i,\\
			B_j,  &\text{ if } c_{ij}=0, \\
			B_jB_i-vB_iB_j,  & \text{ if }c_{ij}=-1, \\
			{[}2]^{-1} \big(B_jB_i^{2} -v[2] B_i B_jB_i +v^2 B_i^{2} B_j \big) + B_j\K_i,  & \text{ if }c_{ij}=-2,
		\end{cases}
		\]
		for $\mu\in \Z\I$ and $j\in \I$.
		
		(2) For $i\in \I$ such that $c_{i,\tau i}=0$, there exists an automorphism $\TT_i$ of the $\Q(v)$-algebra $\tUi$ such that
		$\TT_i(\K_j) =(-v)^{-c_{ij}-c_{\tau i,j}}\K_{j}\K_i^{-c_{ij}} \K_{\tau i}^{-c_{\tau i,j}}$, and
		\[
		\TT_i(B_j)= \begin{cases}
			B_jB_i -vB_iB_j,  & \text{ if }c_{ij}=-1 \text{ and } c_{\tau i,j}=0,
			\\
			B_jB_{\tau i}-v B_{\tau i}B_j,  & \text{ if } c_{ij}=0 \text{ and }c_{\tau i,j}=-1 ,
			\\
			\big[[B_j,B_i]_v,B_{\tau i}\big]_v - v B_j\bK_{i},  & \text{ if } c_{ij}=-1 \text{ and }c_{\tau i,j}=-1 ,
			\\
			-\bK_{i}^{-1} B_{\tau i},  & \text{ if }j=i,
			\\
			-\bK_{\tau i}^{-1} B_i,  &\text{ if }j=\tau i,
			\\
			B_j, & \text{ otherwise;}
		\end{cases}
		\]
		for $\mu\in \Z\I$ and $j\in \I$.

		Moreover, there exists a homomorphism $\brW \rightarrow \Aut( \tUi)$, $r_i\mapsto \TT_i$, for all $i\in \I_\tau$.
	\end{theorem}

        \begin{proof}
        We refer to \cite{Z23} for a precise reference here. The existence of the automorphism $\TT_i$ is established in \cite[Theorem 3.1]{Z23}. The formulas of $\TT_i$ are obtained by specifying those formulas in \cite[Proposition 3.4 and Theorem 3.7]{Z23}. The braid relations for $\TT_i$ are proved in \cite[Theorem 3.11]{Z23}.
        \end{proof}

        \begin{remark}
        One always has that $\TT_i=\TT_{\tau i}$ for any $i\in \I$. For finite type, this is explained in \cite[Remark 4.8]{WZ23} and the argument for Kac-Moody type is essentially the same. Then we can parametrize $\TT_i$ using the set $\I_\tau$ instead of $\I$ and the action of $r_i$ on $\tUi$ does not depend on the choice of representative for each $\tau$-orbit.
        \end{remark}
	
	\begin{remark}
		Let $\Phi$ be the rescaling automorphism on $\tUi$ such that
		\begin{align*}
			B_i\mapsto B_i,\quad B_{\tau i}\mapsto -v^{-1} B_{\tau i}, \quad\K_i\mapsto -v^{-1}\K_i,\quad\K_{\tau i} \mapsto -v^{-1}\K_{\tau i},
		\end{align*}
		for $\tau i\neq i\in \I_\tau$, and $\Phi$ fixes $B_j,\K_j$ for $\tau j =j$. Then the symmetries $\mathrm{T}_i$ in \cite[\S 5]{LW21a} are given by $\mathrm{T}_i=\Phi \TT_i \Phi^{-1}$.
	\end{remark}

 There exists a $\Q(v)$-algebra anti-involution $\sigma_\imath: \tUi\rightarrow \tUi$ such that
\begin{align}
  \label{eq:sigma}
\sigma_\imath(B_i)=B_{i}, \quad \sigma_\imath(\K_i)= \K_{\tau i},
\quad \forall i\in \I.
\end{align}

By \cite[Theorem 6.7]{WZ23}, we have
\begin{align}
\TT_i^{-1}=\sigma_\imath \TT_i\sigma_{\imath}, \quad\forall i\in\I.
\end{align}

The following was proved in \cite[Theorem 7.13]{WZ23} for $\tUi$ of arbitrary finite type, and we extend it to the quasi-split affine type in our setting.
	\begin{lemma}
		\label{lem:fixB}
		We have $\TT_w (B_i) = B_{w i}$, for $i\in \I$ and $w \in W_\tau$ such that $wi \in \I$.
	\end{lemma}
	
	\begin{proof}
	    By similar arguments as in \cite[Theorem 7.13]{WZ23}, one reduces the proof of this statement to the cases of an arbitrary rank two Satake subdiagram $\mathbb J=\{i,\tau i,j,\tau j\}$ of $\I$. If the rank two diagram $\mathbb J$ is of finite type, this statement has been established {\em loc. cit.}. If the rank two diagram is not of finite type, then by inspection there do not exist $1\neq w\in W_\tau$ and $k\in \mathbb J$ such that $w k \in \mathbb J$. Hence the lemma holds.
	\end{proof} 
 
 	\begin{lemma}
		\label{lem:Tomeij}
		Let $i,j\in \I_{\tau}\setminus\{0\}$. Then we have
		\begin{itemize}
			\item[(1)] $\TT_{\bome_i}\TT_{\bome_j}=\TT_{\bome_j}\TT_{\bome_i}$;
			\item[(2)] $\TT_i^{-1} \TT_{\bome_i} \TT_i^{-1}=\TT_{\bome_i}^{-1} \prod_{k\neq i,k\in \I_\tau} \TT_{\bome_k}^{-\ov{c}_{ik}}$.
		\end{itemize}
	\end{lemma}
	
	\begin{proof}
		These formulas follow by repeating the proof of Lemma~\ref{lem:Braid-Lus} (which was for the usual affine braid group) for the relative affine braid group $\brW$.
	\end{proof}
For $i\in\II$, recalling $t_{\bome_i}$ from \eqref{tbome2} and defining
\begin{align*}  
t_{\bome_i'}:=t_{\bome_i}\bs_i, 
\end{align*}
we have that 
\begin{align}  \label{bome2}
    \ell(t_{\bome_i'}) =\ell(t_{\bome_i})-\ell(\bs_i), \qquad \text{and} \qquad \TT_{\bome_{i}'}=\TT_{\bome_i} \TT_i^{-1}. 
\end{align}
	\begin{lemma}
		\label{lem:T-T}
		We have
		$\TT_{\bome_i'}(B_i) =\TT_{\bome_i'}^{-1}(B_i),$ for $i\in \II$.
	\end{lemma}
	
	\begin{proof}
		By Lemma~\ref{lem:Tomeij}, we have $\TT_{\bome_i'}=\TT_{\bome_i'}^{-1} \prod_{j\neq i,j\in \I_\tau} \TT_{\bome_j}^{-\ov{c}_{ij}}$. Then, by Lemma~\ref{lem:fixB}, we have $\TT_{\bome_i'}(B_i)=\TT_{\bome_i'}^{-1} \prod_{j\neq i,j\in \I_\tau} \TT_{\bome_j}^{-\ov{c}_{ij}}(B_i)=\TT_{\bome_i'}^{-1}(B_i)$.
	\end{proof}
 Define a sign function
\[
o(\cdot): \II \longrightarrow \{\pm 1\}
\]
such that $o(i) o(j)=-1$ whenever $c_{ij} <0$ (there are exactly 2 such functions in general).
By inspection of the Satake diagrams in Table \ref{tab:Satakediag}, we see that the values of $o(i)o(\tau i)$ are independent of $i\in \I_0$. Note that the element $ \K_i^{-1} \TT_{\bome_i}^{-1}(\K_i)$ is a scalar multiple of $\K_\delta$, which is independent of $i\in \I_0$. Therefore, the element
\begin{align}
    \label{C}
C:=o(i)o(\tau i)\,\K_i^{-1} \TT_{\bome_i}^{-1}(\K_i)
\end{align}
in $\tUi$ is independent of $i\in \I_0$.

\section{Rank one subalgebras of affine $\imath$quantum groups}
 \label{sec:subalgebras}

We continue to work with the quasi-split affine $\imath$quantum groups $\tUi$ formulated in \S\ref{subsec:qsAffine}. The affine rank one subalgebra of $\tUi$ associated with $i\in \I_0$ is defined to be the subalgebra $\tUi_{[i]}$ with the following generators (see \eqref{bome2}--\eqref{C} for notations): 
\begin{align} \label{def:Uii}
\tUi_{[i]} = \big\langle B_j, \K_j^{\pm1}, \TT_{\bome_j'}(B_j),  C^{\pm1}, \mid j\in\{i,\tau i\} \big\rangle.
\end{align}

In this section we establish an embedding from a suitable affine rank one $\imath$quantum group to $\tUi$, with  $\tUi_{[i]}$ as the image, for each $i\in \I_0$. We then show that the relative braid group actions are compatible with such embeddings.

\subsection{Rank one subalgebras for $c_{i,\tau i}=2$}

Let $i\in \II$ be such that $c_{i,\tau i}=2$, i.e., $\tau i=i$, in this subsection. We consider $\tUi$ together with the affine $\imath$quantum group of split rank one (=q-Onsager algebra) $\tUi(\widehat{\mathfrak{sl}}_2)$; see Remark~\ref{rem:Onsager}. 
To distinguish notations from those for $\tUi$, we shall adopt the dotted notation for $\tUi(\widehat{\mathfrak{sl}}_2)$; that is, $\tUi(\widehat{\mathfrak{sl}}_2) = \langle \dot B_i, \dot \K_i \mid i\in\{0,1\} \rangle$, and  it is acted by the relative braid group operators $\dot \TT_w$.	

Set
	\begin{align}
		\Delta(\omega_i'):= \big \{ \tau_i\alpha_{i_1},\tau_i s_{i_1}(\alpha_{i_2}),\ldots,\tau_is_{i_1}s_{i_2}\cdots s_{i_{k-1}}(\alpha_{i_k}) \big\}.
	\end{align}
	Note that $\Delta(\omega_i')=\{\gamma\in \cR^+ \mid t_{\omega_i'}^{-1}(\gamma) \in -\cR^+\}$ where we recall $\cR^+$ from \eqref{eq:roots+}.
 
\begin{lemma} \label{C:weight}
		We have $\Delta(\omega_i')\subset\{k\delta-\beta\mid k>0,\beta\in \cR_0^+\}$. In particular, the multiplicity of $\alpha_0$ in any root of  $\Delta(\omega_i')$ is positive.
  \end{lemma}
	
\begin{proof} 
		Clearly, any root in $\Delta(\omega_i')$ is a real root. Suppose otherwise that there exists a root $\gamma \in\Delta(\omega_i')$ such that $\gamma\in\{\beta+k\delta\mid \beta\in\cR_0^+,k\geq 0 \}$. Write $\gamma=\sum_{j\in \I_0}c_j \alpha_j +k\delta$ for some $c_j\geq 0,k\geq 0$. Then we have
		\begin{align*}
			t_{\omega_i}^{-1}(\gamma)= \sum_{ j\in \I_0} c_j \alpha_j +(k+c_i)\delta \in \cR^+.
		\end{align*}
		On the other hand, since $t_{\omega_i}^{-1}=s_i t_{\omega_i'}^{-1}$ and $t_{\omega_i'}^{-1}(\gamma)\in -\cR^+$, $t_{\omega_i'}^{-1}(\gamma)$ is a negative root sent to a positive root by $s_i$. Hence, $t_{\omega_i'}^{-1}(\gamma)=-\alpha_i$. However, this implies that $\gamma= t_{\omega_i'}(-\alpha_i)=t_{\omega_i} s_i (-\alpha_i)=\alpha_i-\delta$, which contradicts that $\gamma\in \cR^+$. Thus, this lemma is proved.
\end{proof}

\begin{lemma}\label{lem:weight}
For $\tau i =i\in \II$, we have
\begin{align*}
		&\Big[\big[[ \TT_{\omega'_i}^{-1}(B_i) ,B_i]_{v^2} , B_i\big],B_i\Big]_{v^{-2}}-v[2]^2[\TT_{\omega'_i}^{-1}(B_i),B_i] \K_i=0.
	\end{align*}
\end{lemma}
The lengthy and technical proof of Lemma~\ref{lem:weight} can be found in Appendix \ref{app:A1}.  

\begin{proposition}
		\label{prop:rank1iso}
	For $\tau i =i\in \II$, there is a $\Q(v)$-algebra isomorphism  $\aleph_i: \tUi(\widehat{\mathfrak{sl}}_2) \longrightarrow \tUi_{[i]}$, which sends $\dot B_1 \mapsto B_i$, $\dot B_0 \mapsto \TT_{\bome_i'} (B_i)$, $\dot \K_1 \mapsto \K_i$, and $\dot \K_0 \mapsto \TT_{\bome_i'}(\K_i)$.
\end{proposition}
	
    \begin{proof}
    
	We show that the map $\aleph_i$ is an algebra homomorphism. By \cite[Lemma 3.7]{LW21b} (or Lemma~\ref{lem:T-T}), the identity \eqref{eq:embedding} in Lemma~\ref{lem:weight} implies that
	\begin{align} \label{eq:embedding2}
		&\Big[\big[[  B_i,\TT_{\omega'_i}^{-1}(B_i)]_{v^2} , \TT_{\omega'_i}^{-1}(B_i)\big],\TT_{\omega'_i}^{-1}(B_i)\Big]_{v^{-2}}-v[2]^2[B_i,\TT_{\omega'_i}^{-1}(B_i)] \TT_{\omega'_i}^{-1}(\K_i)=0.
	\end{align}
	It follows by the Serre-type relations  \eqref{eq:embedding}, \eqref{eq:embedding2} and the additional relations involving $\K$'s (which are easily verified) that $\aleph_i: \tUi(\widehat{\mathfrak{sl}}_2) \rightarrow \tUi_{[i]}$ is a homomorphism.
 
    
    The surjectivity of $\aleph_i$ is clear since all generators of $\tUi_{[i]}$ are in the image. We prove the injectivity of $\aleph_i$ as follows (similar to \cite{LW21b}). Recall from \eqref{eq:filt1} and \eqref{eq:filter} the natural filtrations on $\tUi$ (and $\tUi(\widehat{\mathfrak{sl}}_2)$, respectively) such that the associated graded are given by
    \begin{align*} 
    \mathrm{gr}\tUi \cong \U^-\otimes \bQ(v)[\K_i^{\pm1} | i\in \I],
    \qquad
    \mathrm{gr}\tUi(\widehat{\mathfrak{sl}}_2)\cong \U(\widehat{\mathfrak{sl}}_2)^-\otimes \bQ(v)[\dot{\K}_1^{\pm1},\dot{\K}_0^{\pm1}].
    \end{align*}
    The map $\aleph_i$ is compatible with these two filtrations and then $\aleph_i$ induces a homomorphism on the associated graded, $\aleph_i^{\text{gr}}: \U(\widehat{\mathfrak{sl}}_2)^-\rightarrow \U ^-$, which sends $F_1 \mapsto F_i, F_0 \mapsto \tcT_{\omega'_i}(F_i)$. This induced homomorphism $\aleph_i^{\text{gr}}$ was well known in the literature and it is injective \cite{Be94, Da15}, and hence $\aleph_i$ is also injective.
	\end{proof}

\begin{remark}
The node $i\in \II$ satisfying $\tau i=i$ (so Proposition \ref{prop:rank1iso} applies) could be arbitrary if $\tau ={\rm Id}$, or is one of the following if $\tau \neq {\rm Id}$ (see Table \ref{tab:Satakediag}):
	\begin{align*}
		i\in \begin{cases} \{r\}, &\text{ in  AIII}_r^{(\tau)}, \\
			\{1,2,\dots, r-2\}, &\text{ in DI}_r^{(\tau)},
			\\
			\{3,4\},&\text{ in EI}_6^{(\tau)}.  \end{cases}
	\end{align*}
\end{remark}

\subsection{Rank one subalgebras for $c_{i,\tau i}=0$}

Denote by $\ov{\U}(\widehat{\mathfrak{sl}}_2)$ the quotient algebra of $\tU(\widehat{\mathfrak{sl}}_2)$ modulo the ideal generated by the central element $\tK_\de-\tK'_\de$; the algebra $\ov{\U}(\widehat{\mathfrak{sl}}_2)$ is a variant of the Drinfeld double of quantum $\widehat{\mathfrak{sl}}_2$. As we shall consider $\ov{\U}(\widehat{\mathfrak{sl}}_2)$ and $\tUi$ together, we would like to distinguish notations. We shall adopt the dotted notation for the generators of $\ov{\U}(\widehat{\mathfrak{sl}}_2)$: $\dot F_a, \dot E_a, \dot K_a, \dot K_a'$, for $a\in \{0,1\}$, and denote by $\dot\TT_a$ and $\dot\TT_w$ Lusztig's braid group symmetries of $\ov{\U}(\widehat{\mathfrak{sl}}_2)$.  

We set
	\begin{align}
		\Delta(\bome_i'):=\big\{ \tau_i \balpha_{i_1},\tau_i \bs_{i_1}(\balpha_{i_2}),\ldots,\tau_i \bs_{i_1}\bs_{i_2}\cdots \bs_{i_{k-1}}(\balpha_{i_k}) \big\}.
	\end{align}
	Note that $\Delta(\bome_i')=\{\ov{\gamma}\mid  \gamma\in \cR^+, t_{\bome_i'}^{-1}(\ov{\gamma})\in - \ov{\cR}^+\}\subset \ov{\cR}^+$ (recall $\balpha_j,\ov{\cR}^+$ from \S \ref{sub:Weyl}).
	Since the vertex $0$ is fixed by $\tau$, we have $\balpha_0=\alpha_0$ and $\ov{\delta}=\delta$.
 
    \begin{lemma}
		\label{C:qsweight}
		We have $\Delta(\bome_i')\subset\{k\delta-\ov{\beta} \mid k>0,\beta\in \cR_0^+\}$. In particular, the multiplicity of $\alpha_0$ in any root in $\Delta(\omega_i')$ is positive.
    \end{lemma}

\begin{proof}	
 Suppose that there exists an element $\ov{\gamma} \in\Delta(\bome_i')$ such that $\gamma\in\{\beta+k\delta| \beta\in\cR_0^+,k\geq 0 \}$. Write $\gamma=\sum_{j\in \I_0}c_j \alpha_j +k\delta$, for some $c_j\geq 0,k\geq 0$, and then $\ov{\gamma}=\sum_{j\in \I_0}c_j \balpha_j +k\delta$. We have
	$t_{\bome_i}^{-1}(\ov{\gamma})= \sum_{ j\in \I_0} c_j \balpha_j +(k+c_i+c_{\tau i})\delta \in \ov{\cR}^+.$
On the other hand, since $t_{\bome_i}^{-1}=\bs_i t_{\bome_i'}^{-1}$ and $t_{\bome_i'}^{-1}(\ov{\gamma}) \in -\ov{\cR}^+$, $t_{\bome_i'}^{-1}(\ov{\gamma})$ is a negative root in $\ov{\cR}$ which is sent to a positive root by $\bs_i$. Hence, $t_{\bome_i'}^{-1}(\ov{\gamma})=-l\balpha_i$ for some $l>0$ (recall that $\ov{\cR}$ is not necessary reduced). However, this implies $\ov{\gamma}= lt_{\bome_i'}(-\balpha_i)=lt_{\bome_i} \bs_i (-\balpha_i)=l(\balpha_i-\delta)$, which contradicts that $\ov{\gamma}\in \ov{\cR}^+$. Hence, Lemma~\ref{C:qsweight} follows.
\end{proof}

\begin{lemma}\label{lem:qsweight}
For $i\in\II$ with $c_{i,\tau i}=0$, we have
\begin{align}
		\label{eq:qsembedding1}
		\Big[\big[[ \TT_{\bome'_i}^{-1}(B_i) ,B_i]_{v^2}, B_i\big],B_i\Big]_{v^{-2}} &=0,\\
		[\TT_{\bome'_i}^{-1}(B_i),B_{\tau i}] &=0,
  \label{eq:qsembedding1a}\\
  [B_i,B_{\tau i}] &=\frac{\K_{\tau i}-\K_i}{v-v^{-1}}. 
		\label{eq:qsembedding2}
	\end{align}
\end{lemma}

The somewhat lengthy and technical proof of Lemma~\ref{lem:qsweight} can be found in Appendix \ref{app:A2}.

Recall the subalgebra $\tUi_{[i]}$ of $\tUi$ from \eqref{def:Uii}. The following proposition can be viewed as a version of \cite[Proposition 3.8]{Be94}. 
\begin{proposition}
		\label{prop:rank1isoQG}
	For $i\in\II$ with $c_{i,\tau i}=0$,  there is a $\Q(v)$-algebra isomorphism  $\aleph_i: \ov{\U}(\widehat{\mathfrak{sl}}_2) \longrightarrow \tUi_{[i]}$, which sends $\dot F_1 \mapsto B_i, \dot F_0 \mapsto \TT_{\bome_i'} (B_i), \dot E_1\mapsto B_{\tau i}, \dot E_0\mapsto \TT_{\bome_i'}(B_{\tau i}), \dot K_1 \mapsto \K_i, \dot K'_1\mapsto \K_{\tau i}, \dot K_0 \mapsto \TT_{\bome_i'}(\K_i), \dot K'_0\mapsto \TT_{\bome_i'}( \K_{\tau i})$. 
\end{proposition}
	
\begin{proof}	
Assuming that $\aleph_i: \ov{\U}(\widehat{\mathfrak{sl}}_2) \rightarrow \tUi_{[i]}$ is a homomorphism, then the injectivity of $\aleph_i$ is proved as in \cite[Proposition 3.8]{Be94}, and the surjectivity of $\aleph_i$ is clear. Hence it suffices to show that $\aleph_i$ is a homomorphism. 

By definition, when $c_{i,\tau i}=0$, we have $t_{\bome_i'}=t_{\bome_i} \bs_i=t_{\omega_i }t_{\omega_{\tau i}} s_i s_{\tau i}$. Recall identities \eqref{eq:qsembedding1}--\eqref{eq:qsembedding2} in Lemma~\ref{lem:qsweight}. 
By Lemma~\ref{lem:T-T}, applying $\TT_{\bome'_i}^{-1}$ to \eqref{eq:qsembedding1}--\eqref{eq:qsembedding2} gives us the following identities:
	\begin{align}
		\label{eq:qsembedding1'}
		\Big[\big[[B_i, \TT_{\bome'_i}^{-1}(B_i)]_{v^2},  \TT_{\bome'_i}^{-1}(B_i)\big], \TT_{\bome'_i}^{-1}(B_i)\Big]_{v^{-2}} &=0,
		\\
		[B_i,\TT_{\bome'_i}^{-1}(B_{\tau i})] &=0,
  \\
  [\TT_{\bome'_i}^{-1}
  (B_i),\TT_{\bome'_i}^{-1}(B_{\tau i})]
		&= \frac{\TT_{\bome'_i}^{-1}(\K_{\tau i})
			-\TT_{\bome'_i}^{-1}(\K_{i})}{v-v^{-1}}.
		\label{eq:qsembedding2'}
	\end{align}
It follows by the identities \eqref{eq:qsembedding1}--\eqref{eq:qsembedding2'} and the additional easy relations involving the $\K$'s that $\aleph_i: \ov{\U}(\widehat{\mathfrak{sl}}_2) \rightarrow \tUi_{[i]}$ is a homomorphism. 
\end{proof}

\begin{remark}
    The node $i$ such that $c_{i,\tau i}=0$ (and so Proposition \ref{prop:rank1isoQG} applies) is one of the following; see Table \ref{tab:Satakediag}:
	\begin{align*}
	i\in \begin{cases} 
\{1,2,\dots,r-1,r+1,\dots,2r\}, &\text{ in  AIII}_r^{(\tau)},
   \\
\{r-1,r\}, &\text{ in  DI}_r^{(\tau)},
			\\
	\{1,2,5,6\},&\text{ in  EI}_6^{(\tau)}.
		\end{cases}
	\end{align*}	
\end{remark}

\subsection{Braid group actions and rank one embeddings}

	\begin{theorem}  [\cite{WZ24}]
		\label{rank1T}
	For $i\in\II$, the action of $\TT_i$ on any integrable $\U$-module admits a closed formula which depends only on $B_i, B_{\tau i}$ and $\bK_i$.
	\end{theorem}
	We do not need the precise rank one formulas (other than the existence) of the action of $\TT_i$ on modules in this paper; the formulas are very different, depending on whether $c_{i,\tau i} =2$ or $c_{i,\tau i} =0$. In the latter case, it is reduced to Lusztig's formulas \cite{Lus93}.
	
	Recall the subalgebra $\tUi_{[i]}$ of $\tUi$ from \eqref{def:Uii}. Recall from Propositions~\ref{prop:rank1iso} and \ref{prop:rank1isoQG} the embeddings $\aleph_i$ of the affine rank one subalgebras of 2 different types onto $\tUi_{[i]}$.
	We have the following proposition which generalizes \cite[Corollary in  p.560]{Be94} in case $c_{i,\tau i} =0$. In case $c_{i,\tau i} =2$ (for $\tUi$ of split type), this was stated as \cite[(3.27)]{LW21b} without proof.
	
	\begin{proposition}
		\label{prop:T1Ti}
		For $i\in\II$, we have
		\begin{align}
			\TT_i|_{\tUi_{[i]}} &= \aleph_i \circ {\dot \TT}_1 \circ \aleph_i^{-1},
			\label{eq:T1Ti}			\\
			\TT_{\bome_i}|_{\tUi_{[i]}} &= \aleph_i \circ {\dot \TT}_{\bome_1} \circ \aleph_i^{-1}.
			\label{eq:T1Ti2}
		\end{align}
	\end{proposition}

\begin{proof}
 For any integrable $\U$-module $M$, $\tUi$ acts via the composition of the central reduction $\pi_\bvs: \tUi \rightarrow \Ui_\bvs$ and the embedding $\imath: \Ui_\bvs \rightarrow \U$.
	In all cases, by Theorem~\ref{rank1T}, we have
	\begin{align}  \label{eq:T1TiM}
	   \aleph_i \circ {\dot \TT_1} \circ \aleph_i^{-1} =\TT_i \qquad \text{ when acting on } M.
	\end{align}
Recall the compatibility of the actions of $\TT_i$ on $\tUi$ and on modules, i.e., for any $u\in \tUi$ and $x \in M$,  we have
	\begin{align} \label{eq:Tiux}
	\TT_i (u x) = \TT_i(u) \TT_i(x),
	\qquad
	   \aleph_i {\dot \TT}_1 \aleph_i^{-1} (ux)= \aleph_i {\dot \TT}_1 \aleph_i^{-1} (u) \aleph_i {\dot \TT}_1 \aleph_i^{-1}(x).
	\end{align}
Indeed, for any $x \in M$, we denote $x' =\TT_i(x)$, and thus for any $u\in \tUi_{[i]}$, by applying  \eqref{eq:T1TiM} twice and \eqref{eq:Tiux}, we have
	\begin{align}
	  \aleph_i {\dot \TT}_1 \aleph_i^{-1} (u)\cdot x'
	    &= \aleph_i {\dot \TT}_1 \aleph_i^{-1} (u)\cdot \aleph_i {\dot \TT}_1 \aleph_i^{-1} (x)
	    \notag \\
	    &= \aleph_i {\dot \TT}_1 \aleph_i^{-1} (ux)
	    = \TT_i(ux)
	    = \TT_i(u) x'.
	    \label{eq:T1TiUM}
	\end{align}
	Note that, for a given $\tilde u\in \tUi$, if  $\pi_\bvs(\tilde u)=0$ for each $\bvs$, then $\tilde u=0$. Therefore, (a variant of) \cite[3.5.4]{Lus93} is applicable and we conclude from \eqref{eq:T1TiUM} that $\aleph_i {\dot \TT}_1 \aleph_i^{-1} (u) =\TT_i(u)$, for all $u\in \tUi_{[i]}$. This proves the formula \eqref{eq:T1Ti}.
	
It remains to prove the identity \eqref{eq:T1Ti2}. We proceed case-by-case: $c_{i,\tau i} =2$ or $c_{i,\tau i} =0$.  The case when $c_{i,\tau i} =0$ is the same as in \cite{Be94} and will be skipped. 

We focus now on the remaining case when $\tau i =i$. Recall from Proposition~\ref{prop:rank1iso} that
	$\aleph_i: \tUi(\widehat{\mathfrak{sl}}_2) \rightarrow \tUi_{[i]}$ sends $\dot B_1 \mapsto B_i, \dot B_0 \mapsto \TT_{\bome_i'} (B_i), \dot \K_1 \mapsto \K_i, \dot \K_0 \mapsto \TT_{\bome_i'}(\K_i)$.
		
	Note by Lemma \ref{lem:T-T} that $\TT_{\bome_i'}(B_i) =\TT_{\bome_i'}^{-1}(B_i)$. Applying \eqref{eq:T1Ti}, we have
	\begin{align*}
	    \TT_{\bome_i}^{-1} \big(\TT_{\bome_i'}(B_i)\big) &=\TT_i^{-1}(B_i),
	    \qquad
	    \TT_{\bome_i}^{-1} (B_i)
	    =\TT_i^{-1}\big(\TT_{\bome_i'} (B_i)\big).
	\end{align*}
	On the other hand, recall from \cite[\S2]{LW21b} that $\dot \TT_{\bome_1} =\dag \circ \dot \TT_1$ (where $\dag$ is the diagram involution) acting on $\tUi(\widehat{\mathfrak{sl}}_2)$, and we have
		\begin{align*}
	    \dot \TT_{\bome_1}^{-1} (\dot B_0) &= \dot \TT_1^{-1}(\dot B_1),
	    \qquad
	    \dot \TT_{\bome_1}^{-1} (\dot B_1) = \dot \TT_1^{-1}(\dot B_0).
	\end{align*}
The identity \eqref{eq:T1Ti2} follows by comparing the above four formulas (and the simple formulas for relative braid group actions on the $\K$'s which we skip).
\end{proof}

\begin{remark}
In case $\tau i =i$, the formulas \cite[(2.6)-(2.7), (2.13)]{LW21b} for $\dot \TT_{\bome_1}$ acting on $\tUi(\widehat{\mathfrak{sl}}_2)$ can be transformed via \eqref{eq:T1Ti2} to the following formulas:
\begin{align*}
\TT_i \big( \TT_{\bome_i'}(B_i) \big)
&= \frac{1}{[2]}
\big( \TT_{\bome_i'}(B_i) B_i^{2} -v [2] B_i \TT_{\bome_i'}(B_i) B_i +v^2 B_i^{2}  \TT_{\bome_i'}(B_i) \big) + \TT_{\bome_i'}(B_i) \K_i,
\\
\TT_i^{-1} \big( \TT_{\bome_i'}(B_i) \big)
&= \frac{1}{[2]}
\big(B_i^{2} \TT_{\bome_i'}(B_i) -v [2] B_i \TT_{\bome_i'}(B_i) B_i +v^2 \TT_{\bome_i'}(B_i) B_i^{2} \big) + \TT_{\bome_i'}(B_i) \K_i,
\\
\TT_i \big([\TT_{\bome_i'}(B_i), B_i]_{v_i^{-2}} \big)
&=[B_i, \TT_{\bome_i'}(B_i)]_{v_i^{-2}}.
\end{align*}
Similar rank two identities between real root vectors can be written down for the other type of $i\in \I_0$. 
\end{remark}
	
\section{A Drinfeld type presentation of quasi-split affine $\imath$quantum groups}
 \label{sec:Dpresentation}
	
In this section, we define various real and imaginary root vectors for $\tUi$. Then we formulate a Drinfeld type presentation for the quasi-split affine $\imath$quantum group $\tUi$ using these root vectors. 

\subsection{Root vectors in higher ranks}
	\label{subsec:root vectors}
	Define a sign function
	\[
	o(\cdot): \II \longrightarrow \{\pm 1\}
	\]
	such that $o(j)=-o(i)$ whenever $c_{ij} <0$. (There are clearly exactly 2 such functions.) We define uniformly the following elements in $\tUi$ (called {\em real root vectors}), for $i\in \II$, $k\in \Z$: 
	\begin{align}
 \label{re root}
		B_{i,k} &= o(i)^{-k} \TT_{\bome_i}^{-k} (B_i).
	\end{align}
	In particular, we have $B_{i,0}=B_i$.

 Next we define case-by-case the imaginary root vectors $\TH_{i,n}$, for $i\in \II$, $n\ge 1$,  depending on whether $c_{i,\tau i}=2$ or $0$.

 \subsubsection{The case when $c_{i,\tau i}=2$}
 
 For $n\in \Z$ and $i\in\I_0$ such that $c_{i,\tau i}=2$, we denote
	\begin{equation}\label{GDn}
		D_{i,n }=-[B_{\tau i},B_{i,n }]_{v^{c_{i,\tau i}}}-[B_{i,n+1},B_{\tau i,-1}]_{v^{c_{i,\tau i}}}.
	\end{equation}
Set $\TH_{i,0}=\frac{1}{v-v^{-1}}$. 	Define $\TH_{i,n}$, for $n\ge 1$, recursively as follows:
	\begin{align}\label{GTH2}
		\TH_{i,1} &=  o(i) (-B_{i} \TT_{\bome_i'} (B_i) +v^{2} \TT_{\bome_i'} (B_i) B_{i}),
	\\
		\label{eq:GTH12}
		\TH_{i,2}&=D_{i,0} C \bK_{i}^{-1}-v^{-1}C,
	\\
 \label{GTHn}
		\TH_{i,n}\bK_{i} &=v^{-2}\TH_{i,n-2}\bK_{ i}C +D_{i,n-2}C. 
	\end{align}
We also set  $\TH_{i,n}=0$ for $n<0$. These definitions are inspired by the formulas for root vectors in the split affine type \cite{LW21b} (also cf. \cite{BK20}). 

 \subsubsection{The case when $c_{i,\tau i}=0$}
	For $i\in\I_0$ with $c_{i,\tau i}=0$, we define
	\begin{align}
		\label{def:thetac=0}
		\Theta_{i,0}:=\frac{1}{v-v^{-1}},\qquad
		\Theta_{i,n}:=[B_{i,n},B_{\tau i}]\K_{ \tau i}^{-1}, \text{ for }n>0.
	\end{align}	
	
By Proposition~\ref{prop:T1Ti}, $\Theta_{i,n}=o(i)^{n}[\TT_{\bome_i}^{-n}(B_i),B_{\tau i}]\K_{\tau i}^{-1}$ is identified with the $\aleph_i$-image of $[x_{1,-n}^-,x_{1,0}^+]K_1'^{-1} =\frac{\varphi_{1,-n}}{v-v^{-1}}$; see Proposition~\ref{prop:rank1isoQG} for $\aleph_i$.

 Back to the general setting, we now formulate a key property shared by imaginary root vectors in all types.

	\begin{proposition} \label{invTheta}
 The following identity holds:
 \[
 \TT_{\bome_j}(\Theta_{i,n})=\Theta_{i,n},
 \]
  for all $i, j\in\II$ and $n \ge 1$.
	\end{proposition}
	
	\begin{proof}
		The statement for $j\neq i$ follows by Lemma~\ref{lem:fixB} and the definition of $\Theta_{i,n}$ above in terms of $B_{i,k}$ and $B_{\tau i,k}$, for various $k$.
		
		The statement for $j=i$ is reduced to the affine rank one case by Proposition~\ref{prop:T1Ti}. The statement in all 2 types of affine rank one is known to hold: the case for $c_{i,\tau i}=0$ was proved in \cite{Da93}, the case for $c_{i,\tau i}=2$ was proved in \cite{BK20, LW21b}.
	\end{proof}

\subsection{A Drinfeld type presentation}

  Introduce the following shorthand notation:
\begin{align}
		\label{eq:SS}
	\bS_{i,j}(k_1,k_2|l) : = \Sym_{k_1,k_2}\Big(B_{i,k_1}B_{i,k_2}B_{j,l}-[2]B_{i,k_1} B_{j,l} B_{i,k_2} + B_{j,l} B_{i,k_1}B_{i,k_2}\Big),
\end{align}
for $k_1,k_2,l \in \Z$ and $i\neq j \in \II$.

\begin{definition} 
		\label{def:iDRA1}
	Let $\tUiD$ be the $\Q(v)$-algebra generated by the elements $B_{i,l}$, $H_{i,m}$, $\bK_i^{\pm1}$, $C^{\pm1}$, where $i\in\II$, $l\in\Z$ and $m>0$, subject to the following relations \eqref{qsiDR0}--\eqref{qsiDR8}, for $i,j\in \II$, $m,n>0,$ and $ l,k,k_1,k_2\in \Z$:
\begin{align}
	C^{\pm1} \text{ is central,} \quad [\bK_i,\bK_j]&  =  [\bK_i,H_{j,n}]=0,\quad \bK_iB_{j,l}=v^{c_{\tau i,j}-c_{ij}} B_{j,l} \bK_i,
   \label{qsiDR0}
			\\
   \label{qsiDR1}
			 [H_{i,m},H_{j,n}] &=0,
			\\
   \label{qsiDR2}
			[H_{i,m},B_{j,l}] &=\frac{[mc_{ij}]}{m} B_{j,l+m}-\frac{[mc_{\tau i,j}]}{m} B_{j,l-m}C^m,
			\\
   \label{qsiDR4}
		[B_{i,k},B_{\tau i,l}]
			&= \bK_{\tau i} C^l \TH_{i,k-l}- \bK_{i} C^k \Theta_{\tau i,l-k},
			 \text{ if }c_{i,\tau i}=0,
			\\
\notag [B_{i,k},B_{i,l+1}]_{v^{-2}}-&v^{-2}[B_{i,k+1},B_{i,l}]_{v^{2}}
		= v^{-2} \Theta_{i,l-k+1}\bK_{i} C^k -  v^{-4} \Theta_{i,l-k-1}\bK_{i} C^{k+1}
			\\
		+v^{-2} & \Theta_{i,k-l+1} \bK_{i} C^l -v^{-4} \Theta_{i,k-l-1}\bK_{i} C^{l+1},  \text{ if } \tau i =i,
		\label{qsiDR5}
  \\
			%
   \label{qsiDR3}
			[B_{i,k},B_{j,l+1}]_{v_i^{-c_{ij}}} -&v_i^{-c_{ij}}[B_{i,k+1},B_{j,l}]_{v_i^{c_{ij}}} =0, \qquad \text{ if } j\neq \tau i,
\end{align}
and the Serre relations
\begin{align}
	&\bS_{i,j}(k_1,k_2|l)
	= 0, \qquad\qquad\qquad \text{ if }c_{i,j}=-1, j \neq \tau i\neq i, \label{qsiDR9}
  \\
	&\bS_{i,j}(k_1,k_2|l)  
		 =-\sum_{p=0}^{\lfloor\frac{k_2-k_1-1}{2}\rfloor} v^{2p}  [2] [\Theta_{i,k_2-k_1-2p-1},B_{j,l-1}]_{v^{-2}} \bK_{i}C^{k_1+p+1} 
		\notag	\\
	&\qquad\qquad\qquad  -\sum_{p=1}^{\lfloor\frac{k_2-k_1}{2}\rfloor}v^{2p-1}[2] [B_{j,l},\Theta_{i,k_2-k_1-2p}]_{v^{-2}}\bK_{i}C^{k_1+p}
			\notag
			\\
	&\qquad\qquad\qquad  -[B_{j,l}, \Theta_{i,k_2-k_1}]_{v^{-2}}   \bK_{i}C^{k_1}+\{k_1\leftrightarrow k_2\}, \text{ if }c_{ij}=-1 \text{ and }\tau i=i.
			   \label{qsiDR8}
\end{align}
\end{definition}
Here and below, $\{k_1\leftrightarrow k_2\}$ stands for repeating the previous summand with $k_1,k_2$ swapped.

In the definition above, $H_{i,m}$ are related to $\Theta_{i,m}$ by the following equation for generating functions in $u$: 
	\begin{align}
 \label{exp}
		1+ \sum_{m\geq 1} (v-v^{-1})\Theta_{i,m} u^m  = \exp\Big( (v-v^{-1}) \sum_{m\geq 1} H_{i,m} u^m \Big).
	\end{align}

\begin{proposition}
		\label{prop:qsiDR7}
		Assume $c_{ij}=0$, for $i,j \in \II$ such that $\tau i\neq j$. Then the following relation holds in $\tUiD$:
		$[B_{i,k},B_{j,l}]=0,$ for all $k,l\in\Z$.
\end{proposition}  

\begin{proof}
    Follows from the relation \eqref{qsiDR3}. 
\end{proof}
	
\begin{lemma}  \label{lem:aut}
For each $i\in \II$, there exists an algebra automorphism $\t_i$ on $\tUiD$ such that
\begin{align*}
\t_i (B_{i,r}) &=o(i) B_{i,r -1}, \quad \t_i (B_{\tau i,r}) =o(\tau i)B_{\tau i,r -1},\quad \t_i (H_{j,m}) =H_{j,m}, \quad \t_i(C) =C,
\\
\t_i (\K_i) &=o(i)o(\tau i) \K_i C^{-1}, \quad \t_i (\K_{\tau i}) =o(i)o(\tau i) \K_{\tau i} C^{-1}, \quad \t_i(B_{j,r})=B_{j,r},
\end{align*}
 for all $r\in \Z, m\ge 1$, and $j\not \in \{i,\tau i\}$. Moreover, $\t_i=\t_{\tau i}$ and $\t_i \t_k =\t_k\t_i$ for all $i,k \in \II$.
\end{lemma} 

\begin{proof}
    Follows by inspection of the defining relations for $\tUiD$. 
\end{proof}

We can now formulate the main result of this paper. Recall the root vectors $B_{i,k}, \TH_{i,m}$ in $\tUi$ from \S\ref{subsec:root vectors}. We hope that the following theorem convinces the reader the convenience (instead of confusion) of using identical notations for these root vectors and the generators in $\tUiD$. 

\begin{theorem}
		\label{thm:Dr}
There is a $\Q(v)$-algebra isomorphism ${\Phi}: {}\tUiD \longrightarrow\tUi$, which sends
	\begin{align}
		\label{def:Phi}
	B_{i,k}\mapsto B_{i,k}, \quad \Theta_{i,m} \mapsto \Theta_{i,m},
			\quad
	\K_i\mapsto \K_i, \quad C\mapsto o(i)o(\tau i)\K_i^{-1}\TT_{\bome_i}^{-1}(\K_i),
	\end{align}
 for $i\in \II, m\ge 1,$ and $k\in \Z$.
\end{theorem}

\begin{proof}
We assume that $\Phi\colon \tUiD \rightarrow\tUi$ is an algebra homomorphism for now, postponing its proof until Section \ref{sec:relation}. 

The injectivity of $\Phi$ is proved by the same argument as in \cite[Proof of Theorem~ 3.13]{LW21b} by passing to the associated graded, where $\Phi^{\text{gr}}$ becomes a well-known isomorphism \cite{Be94, Da15} for the Drinfeld presentation of half the affine quantum group. We skip the detail.

We show that $\Phi: \tUiD \rightarrow\tUi$ is surjective. All generators for $\tUi$ except $B_0$ are clearly in the image of $\Phi$, and so it suffices to show that $B_0\in \text{Im} (\Phi)$. We shall adapt the arguments in the proof of \cite[Theorem~12.11]{Da12}; see also \cite[Proof of Theorem~ 3.13]{LW21b}.

The automorphisms $\t_i \in \Aut (\tUiD)$  and $\TT_{\bome_i} \in \Aut (\tUi)$, for $i\in \I$, satisfy
\[
\Phi \circ \t_i =\TT_{\bome_i} \circ \Phi, \qquad \text{ for } i\in \I.
\]
It follows by Lemma \ref{lem:aut} that $\text{Im}(\Phi)$ is $\TT_{\bome_i}$-stable, for each $i \in \I$.

Recall $\theta$ is the highest root in $\cR^+_0$. For any $\alpha \in \cR$, $\balpha=\frac{\alpha+\tau \alpha}{2}$ is a root in the relative root system. Note that ${\boldsymbol\theta}=\theta,\balpha_0=\alpha_0,{\boldsymbol \de}=\de$. Write $\theta=\sum_{i\in \I_\tau\setminus\{0\}}r_i\balpha_i.$  
We can choose ${\bf i} \in \II$ with $\tau {\bf i}={\bf i}$ such that either $r_{\bf i}=1$ or $r_{\bf i}=2, c_{0,{\bf i}}=-1$ (and so $\alpha_0+\alpha_{\bf i}$ is a root). Then $\theta_{\bf i}:=\theta-(r_{\bf i}-1)\alpha_{\bf i}$ is a root.

Write $t_{\bome_{\bf i}} =\sigma_{\bf i} \bs_{i_1} \ldots \bs_{i_N}$ with $\sigma_{\bf i} \in \Omega$. Then $t_{\bome_{\bf i}}(\theta_{\bf i}) =\theta_{\bf i}-\de\in -\cR^+$, and so there exists $p$ such that $\bs_{i_N} \ldots \bs_{i_{p+1}} (\balpha_{i_p}) =\theta_{\bf i}$. Set $y:=\TT_{i_N}^{-1} \ldots \TT_{i_{p+1}}^{-1} (B_{i_p}) \in  \text{Im} (\Phi)$ and $\beta:=\wt (y)$; note that ${\boldsymbol \beta} =\theta_{\bf i}$.
Consider $\TT_{\bome_{\bf i}} (y) =  \sigma_{\bf i} \TT_{i_1} \ldots \TT_{i_{p-1}} (B_{i_p}) \in \text{Im} (\Phi)$, since $\text{Im}(\Phi)$ is $\TT_{\bome_{\bf i}}$-stable. Note that $\sigma_{\bf i} s_{i_1} = s_0 \sigma_{\bf i}$ (see \cite[Lemma 3.1]{Be94}), and so $\TT_{s_0t_{\bome_{\bf i}}} =\TT_0^{-1} \TT_{\bome_{\bf i}}$.

If $r_{\bf i}=1$, then $s_0t_{\bome_{\bf i}}({\boldsymbol \beta})=\alpha_0$, and it follows that $\wt \big(\TT_0^{-1} \TT_{\bome_{\bf i}} (y)\big)=\alpha_0$. By Lemma~\ref{lem:fixB}, we have that $\TT_0^{-1} \TT_{\bome_{\bf i}} (y) =B_0$, and hence $\TT_{\bome_{\bf i}} (y) =\TT_0(B_0) =\K_{0} B_0$. 
Therefore, we have $\K_{0} B_0\in \text{Im} (\Phi)$ and thus $B_0\in \text{Im} (\Phi)$. 
If $r_{\bf i}=2$, then $ \TT_{\bome_{\bf i}} (y)$ has weight $-(\alpha_0+\alpha_{\bf i})$. Since $l^\circ(\bs_{\bf i} t_{\bome_{\bf i}})=l^\circ(t_{\bome_{\bf i}})+1$, we have $ \TT_0^{-1}\TT_{\bf i}\TT_{\bome_{\bf i}} =\TT_{s_0 s_{\bf i} t_{\bome_{\bf i}}}$. Now  $ \TT_0^{-1}\TT_{\bf i}\TT_{\bome_{\bf i}}(y)$ has weight $\alpha_0$, and by Lemma~\ref{lem:fixB}, $ \TT_0^{-1}\TT_{\bf i}\TT_{\bome_{\bf i}}(y) =B_0$. Therefore, we have $\TT_0(B_0)=\K_{0} B_0\in \text{Im} (\Phi)$ and thus $B_0\in \text{Im} (\Phi)$. 

This proves the surjectivity of $\Phi$, and completes the proof of the theorem. 
\end{proof}

\subsection{Drinfeld presentation via generating functions}

Introduce the generating functions
\begin{align}\label{eq:Genfun}
\begin{cases}
\bB_{i}(z) =\sum_{k\in\Z} B_{i,k}z^{k},
\\
 \bTH_{i}(z)  =1+ \sum_{m > 0}(v-v^{-1})\Theta_{i,m}z^{m},
\\
\bH_i(u)=\sum_{m\geq 1} H_{i,m} u^m,
\\
\bDel(z)=\sum_{k\in\Z}  C^k z^k.
\end{cases}
\end{align}
The equation \eqref{exp} can be reformulated in terms of generating functions as 
\begin{align}
    \label{expz}
    \bTH_{i}(z)  = \exp \big((v-v^{-1}) \bH_i(z)\big).
\end{align}
Introduce the following notation
\begin{align*}
& \bS_{i,j}(w_1,w_2|z)
\\
&=\Sym_{w_1,w_2} \big( \bB_{j}(z)\bB_{i}(w_{1})\bB_{i}(w_2) -[2] \bB_{i}(w_{1})\bB_{j}(z)\bB_{i}(w_{2}) +\bB_{i}(w_1)\bB_{i}(w_{2})\bB_{j}(z) \big).
\end{align*}
We can rewrite the defining relations for $\tUiD$ via generating functions in \eqref{eq:Genfun}, and hence reformulate Theorem~\ref{thm:Dr} as follows.

\begin{theorem} \label{thm:DrGF}
$\tUi$ is generated by the elements $B_{i,l}$, $\Theta_{i,m}$, $\bK_i^{\pm1}$, $C^{\pm1}$, where $i\in\II$, $l\in\Z$ and $m>0$, subject to the following relations, for $i,j\in \II$:
\begin{align}
 C \text{ is central, } & \quad [\bK_i,\bK_j] =  [\bK_i,\bTH_{j}(w)] =0,\quad \bK_i\bB_j(w)=v^{c_{\tau i,j}-c_{ij}} \bB_j(w) \bK_i,
\label{Rel1GF} \\
\bTH_i(z) \bTH_j(w) &=\bTH_j(w) \bTH_i(z),
\\
 \bB_j(w)  \bTH_i(z)
 &= \frac{1 -v^{c_{ij}}zw^{-1}}{1 -v^{-c_{ij}}zw^{-1}} \cdot \frac{1 -v^{-c_{\tau i,j}}zw C}{1 -v^{c_{\tau i,j}} zw C}
\bTH_i(z) \bB_j(w), 
\\
[\bB_i(z),\bB_{\tau i}(w)] &=\frac{\bDel(zw)}{v-v^{-1}} \big(\K_{\tau i}\bTH_i(z)-\K_i\bTH_{\tau i}(w)\big),
\qquad \text{ if } c_{i,\tau i}=0,
\\
\label{BiBiGF}
(v^2z-w) & \bB_i(z) \bB_i(w) +(v^{2}w-z) \bB_i(w) \bB_i(z)
\\
&= v^{-2} \frac{\bDel(zw)}{v-v^{-1}} \big( (v^2z-w)\bTH_i(w) +(v^2w-z)\bTH_i(z) \big)\K_{i},
 \text{ if } i=\tau i,
\\
(v^{c_{ij}}z -w) & \bB_i(z) \bB_j(w)  +(v^{c_{ij}}w-z) \bB_j(w) \bB_i(z)=0, \qquad \text{ if }j\neq \tau i,
\label{BB1}
\end{align}
and the Serre relations
\begin{align}
    \bS_{i,j}(w_1,w_2|z) &=0,\qquad \text{ if } c_{ij}=-1,j\neq \tau i\neq i,
    \\
\bS_{i,j}(w_1,w_2|z)
&= -\frac{\bDel(w_1w_2)}{v-v^{-1}}  \Sym_{w_1,w_2} \left(\frac{[2] z w_1^{-1} }{1 -v^{2}w_2w_1^{-1}}[\bTH_i(w_2),\bB_j(z)]_{v^{-2}} \right.
\notag \\ 
&\qquad\qquad\qquad\qquad\qquad
+ \left.\frac{1 +w_2w_1^{-1}}{1 -v^{2}w_2w_1^{-1}}[\bB_j(z),\bTH_i(w_2)]_{v^{-2}} \right)\K_i, 
\label{SSGF2} \\
\notag
&\qquad\qquad\qquad\qquad\qquad\qquad\qquad\qquad\qquad \text{ if }c_{ij}=-1, i=\tau i.
\end{align}
\end{theorem}

\begin{lemma}   \label{lem:HH}
The following identities are equivalent, for $i,j \in \II$:
\begin{align}
    [H_{i,m},H_{j,n}] &=0,\quad \forall  m,n\geq1;
  \\
  \bH_i(z) \bH_j(w) &=\bH_j(w) \bH_i(z);
  \\
  [\TH_{i,m},\TH_{j,n}] &=0,\quad \forall m,n\geq1;
    \\ 
  \bTH_i(z) \bTH_j(w) &=\bTH_j(w) \bTH_i(z).
\end{align}
\end{lemma}

\begin{proof}
    Follows by \eqref{eq:Genfun}. 
\end{proof}

\begin{lemma}  \label{lem:HB}
The following identities are equivalent, for $i,j \in \II$:
\begin{align}
[H_{i,m},B_{j,l}] &=\frac{[mc_{ij}]}{m} B_{j,l+m}-\frac{[mc_{\tau i,j}]}{m} B_{j,l-m}C^m, \; \forall m\ge 1, k\in\Z; 
\label{HBij} \\
[\bH_i(z), \bB_j(w)] &= \frac1{v-v^{-1}} \ln\left( \frac{1 -v^{c_{ij}}zw^{-1}}{1 -v^{-c_{ij}}zw^{-1}} \cdot \frac{1 -v^{-c_{\tau i,j}}zw C}{1 -v^{c_{\tau i,j}} zw C}
  \right) \bB_j(w);
\label{HBzw}   \\
[\TH_{i,m}, B_{j,k}] &+v^{c_{i,j}-c_{\tau i,j}}[\TH_{i,m-2},B_{j,k}]_{v^{2(c_{\tau i,j}-c_{i,j})}}C
-v^{c_{i,j}}[\TH_{i,m-1}, B_{j,k+1}]_{v^{-2c_{ i,j}}} 
 \notag \\
 &- v^{-c_{\tau i, j}}[\TH_{i,m-1},B_{j,k-1}]_{v^{2c_{\tau i,j}}}C
			=0,
\; \forall m\ge 1, k\in\Z;
\label{THBij}  \\
\bB_j(w)  \bTH_i(z) &= 
  \frac{1 -v^{c_{ij}}zw^{-1}}{1 -v^{-c_{ij}}zw^{-1}} \cdot \frac{1 -v^{-c_{\tau i,j}}zw C}{1 -v^{c_{\tau i,j}} zw C}
  \bTH_i(z) \bB_j(w).
\label{THBzw}
  \end{align}
\end{lemma}

 \begin{proof}
 The same arguments for \cite[Proposition 2.8]{LW21b} carry over, which we outline here. It is straightforward to convert \eqref{HBij} to its generating function form \eqref{HBzw} and to recover \eqref{HBij} from \eqref{HBzw} by comparing the coefficient of $z^mw^l$. In the same way one proves the equivalence between \eqref{THBij} and its generating function form \eqref{THBzw}. One then shows that \eqref{HBzw} and \eqref{THBzw} are equivalent by taking integration/derivation and using \eqref{expz}. 
 \end{proof}

Lemmas~\ref{lem:HH} and \ref{lem:HB} establish the equivalence of some relations in different forms in Theorems \ref{thm:Dr} and \ref{thm:DrGF}. The equivalence of other relations in Theorems \ref{thm:Dr} and \ref{thm:DrGF} can be similarly established as well. This implies the equivalence between Theorems \ref{thm:Dr} and \ref{thm:DrGF} and we will skip the details. An interested reader can look up the proof of an analogous equivalence in \cite[Theorem 5.1]{LW21b}. 

\subsection{A variant of the Drinfeld presentation}	
	
When $i=\tau i$, a variant of imaginary root vectors denoted by $\acute{\Theta}_{i,m}$ in $\tUi$ was given in \cite{LW21b}, which is compatible with the version for $\Ui$ \cite{BK20} in equal parameter (i.e., $\vs_i=\vs_{0}$ for any $i\in \II$):
\begin{align}
\acute{\Theta}_{i,m} :=\begin{cases}
-\sum\limits_{p=0}^{m-1} [B_{i,p},B_{i,m-2-p}]_{v^2} \K_0 & \text{ if }m\geq1,
\\
\frac{1}{v-v^{-1}} & \text{ if }m=0.
\end{cases}
\end{align}
 The alternative Drinfeld type presentation using these imaginary root vectors for $\tUi$ in the split affine type was given in \cite[Theorem 5.4]{LW21b}. Define 
$$\acute{\bTH}_i(z):=1+ \sum_{m > 0}(v_i-v_i^{-1})\acute{\Theta}_{i,m}z^{m}.$$
By \cite[Lemma 2.9]{LW21b}, we have 
\begin{align}
\label{Theta-transform}
\bTH_i(z) =  \frac{1- Cz^2}{1-v^{-2}C z^2}  \acute{\bTH}_i(z).
\end{align}

Replacing $\Theta$ by $\acute{\Theta}$ in computation (but not in notation below), we formulate an alternative Drinfeld type presentation for $\tUi$ as follows (compare Definition \ref{def:iDRA1} and Theorem \ref{thm:Dr}). 

 \begin{theorem}
 \label{thm:Drvariant}
 $\tUi$ is isomorphic to the $\Q(v)$-algebra generated by $B_{i,k}$, $H_{i,m}$, $\bK_i^{\pm1}$, $C^{\pm1}$, where $i\in \II$, $k\in\Z$ and $m \in \Z_{\ge 1}$, subject to the relations \eqref{qsiDR0}-\eqref{qsiDR4}, \eqref{qsiDR3}--\eqref{qsiDR9} and the following two relations \eqref{qsiDR5BK}--\eqref{qsiDR8BK} (in place of \eqref{qsiDR5} and \eqref{qsiDR8}):
\begin{align}\notag
	[B_{i,k}, B_{i,l+1}]_{v^{-2}} & -v^{-2} [B_{i,k+1},B_{i,l}]_{v^{2}}
\\\notag
&= v^{-2} {\Theta}_{i,l-k+1} C^k \K_i-v^{-2} {\Theta}_{i,l-k-1} C^{k+1} \K_i 
\\
&\quad +v^{-2} {\Theta}_{i,k-l+1} C^l \K_i-v^{-2} {\Theta}_{i,k-l-1} C^{l+1} \K_i,
\qquad \text{ if } i=\tau i,	\label{qsiDR5BK} 
\\
\mathbb{S}_{i,j}(k_1,k_2| l)
	& =- \textstyle  C^{k_1} \K_i \Big(\sum_{p\geq0} 
     (v^{-2p-1} +v^{2p+1}) [ {\Theta}_{i,k_2-k_1-2p-1},B_{j,l-1}]_{v^{-2}}C^{p+1} \notag
        \\\notag
& \textstyle
\quad +\sum_{p\geq 1} (v^{-2p} +v^{2p} )  [B_{j,l}, {\Theta} _{i,k_2-k_1-2p}]_{v^{-2}} C^{p}
+[B_{j,l}, {\Theta} _{i,k_2-k_1}]_{v^{-2}}    \Big)
  \notag \\
 &\quad +\{k_1 \leftrightarrow k_2\},
\qquad \text{ if }c_{ij}=-1, j\neq \tau i\neq i,
			\label{qsiDR8BK}
		\end{align}
for $k_1, k_2, k, l \in\Z$, and $i,j\in \II$.
 \end{theorem}

The generating function form of \eqref{qsiDR5BK}--\eqref{qsiDR8BK} can be readily obtained  as follows by substituting \eqref{Theta-transform} into  \eqref{BiBiGF} and \eqref{SSGF2}: 
\begin{align}
(v^2z-w)& \bB_i(z) \bB_i(w) +(v^{2}w-z) \bB_i(w) \bB_i(z)
\\\notag
&=  \frac{\bDel(zw)}{v-v^{-1}} \Big( (v^2z-w)\frac{1- Cw^2}{v^2-C w^2}\bTH_i(w) +(v^2w-z)\frac{1- Cz^2}{v^2-C z^2}\bTH_i(z) \Big)\K_{i},
\\
\bS_{i,j}(w_1,w_2|z)
&= - \frac{1-C w_2^2}{1- v^{-2}Cw_2^2} \frac{\bDel(w_1w_2)}{v-v^{-1}} \times \\
 \Sym_{w_1,w_2} &
\left(\frac{[2] z w_1^{-1} }{1 -v^{2}w_2w_1^{-1}}[\bTH_i(w_2),\bB_j(z)]_{v^{-2}} 
 + \frac{1 +w_2w_1^{-1}}{1 -v^{2}w_2w_1^{-1}}[\bB_j(z),\bTH_i(w_2)]_{v^{-2}} \right)\K_i.
 \notag
\end{align}

	\section{Verification of relations among root vectors}
	\label{sec:relation}
	
	In this section, we prove that ${\Phi}: \tUiD \rightarrow\tUi$ defined by \eqref{def:Phi} is a homomorphism by verifying the defining relations \eqref{qsiDR0}--\eqref{qsiDR8} of $\tUiD$ for the images of the generators under $\Phi$. The relation \eqref{qsiDR0} is clear. 

\subsection{Rank one relations: \eqref{qsiDR1}--\eqref{qsiDR2} for $j=i$ and \eqref{qsiDR4}--\eqref{qsiDR5}}
	
	\begin{proposition}
		\label{prop:iDR31}
		All rank one relations, that is, \eqref{qsiDR1}--\eqref{qsiDR2} for $j=i$ and \eqref{qsiDR4}--\eqref{qsiDR5}, hold in $\tUi$.
	\end{proposition}
	
	\begin{proof} 
	The current relations in the q-Onsager algebra $\tUi(\widehat{\mathfrak{sl}}_2)$ are obtained in \cite[Definition 2.15, Theorem~2.16]{LW21b}. The current relations in the affine rank one quantum group $\tU(\widehat{\mathfrak{sl}}_2)$ are obtained in \cite{Da93}. Then all the rank one relations follow by transporting the corresponding (rank one) relations in $\tUi(\widehat{\mathfrak{sl}}_2)$ and  $\tU(\widehat{\mathfrak{sl}}_2)$ via Propositions \ref{prop:rank1iso} and \ref{prop:rank1isoQG}, respectively.
	\end{proof}

	\subsection{Relation \eqref{qsiDR3}}

 We first establish a lemma which is similar to \cite[Lemma 3.3]{Be94} in the quantum group setting. 
	\begin{lemma} 
		\label{lem:TXij}
		For $j\neq \tau i\in\II$ such that $c_{ij}=-1$, denote
		\begin{align*}
			X_{ij}:= B_jB_i-vB_iB_j.
		\end{align*}
		Then we have
		$\TT_{\bome_i}(X_{ji})=\TT_{\bome_j}(X_{ij}).$
	\end{lemma}
	
	\begin{proof}		
		According to Theorem~\ref{thm:Ti}, if $\ov{c}_{ji}=-1$, we have
		\begin{align*}
			\TT_j^{-1}(B_i)=B_jB_i-vB_iB_j=X_{ij},
			\qquad
			\TT_j(B_i)=B_iB_j-vB_jB_i=X_{ji}.
		\end{align*}
		Observe that, in this case, we always have either $\ov{c}_{ij}=-1$ or $\ov{c}_{ji}=-1$. Without loss of generality, assume that $\ov{c}_{ji}=-1$.
		Then we have
		\begin{align*}
			\TT_{\bome_j}(X_{ij})&=\TT_{\bome_j} \TT_j^{-1}(B_i)= \TT_j\TT_{\bome_j}^{-1}\TT_{\bome_i}(B_i)
			=\TT_{\bome_i}\TT_j(B_i)=\TT_{\bome_i}(X_{ji}).
		\end{align*}
  The lemma is proved. 
	\end{proof}
	
	Now we are ready to establish the relation \eqref{qsiDR3}.
	\begin{proposition}
		\label{prop:iDR3a}
		We have $[B_{i,k}, B_{j,l+1}]_{v^{-c_{ij}}}  -v^{-c_{ij}} [B_{i,k+1}, B_{j,l}]_{v^{c_{ij}}}=0$, for $j\neq \tau i \in \II$ and $k, l \in \Z.$
	\end{proposition}
	
	\begin{proof}
		For $j=i$, it follows by transporting the corresponding relations in $\tU(\widehat{\mathfrak{sl}}_2)$ by using Proposition \ref{prop:rank1iso}.
		
		It remains to consider the case $i\neq j\neq \tau i$.
		
		Assume first that $c_{ij}=0$. Then the desired identity is equivalent to the identity $[B_{i,k}, B_{j,l}]=0$, for all $k,l\in \Z$. The identity for $k=l=0$, i.e., $[B_{i},B_{j}]=0,$ is the defining relation \eqref{eq:S1} for $\tUi$. The identity for general $k,l$ follows by applying $\TT_{\bome_i}^{-k} \TT_{\bome_j}^{-l}$ to the above identity and using Lemma~\ref{lem:Braid-Lus} and Lemma~\ref{lem:fixB}.
		
		Assume now that $c_{ij}=-1$. Note that
		$v[B_{i,k+1},B_{j}]_{v^{-1}}= - o(i)^{k+1} \TT_{\bome_i}^{-(k+1)}(X_{ij}).$
		By using Lemma \ref{lem:TXij}, we have
		\begin{align*}
			[B_{i,k},B_{j,1}]_{v}&= B_{i,k}B_{j,1}-vB_{j,1}B_{i,k}
			\\
			& = o(i)^k o(j) \TT_{\bome_i}^{-k} \TT_{\bome_j}^{-1}(X_{ji})
			= - o(i)^{k+1} \TT_{\bome_i}^{-k} \TT_{\bome_i}^{-1} (X_{ij})
			\\
			&= - o(i)^{k+1} \TT_{\bome_i}^{-(k+1)}(X_{ij})=v[B_{i,k+1},B_{j}]_{v^{-1}}.
		\end{align*}
		So we have obtained that
		$[B_{i,k},B_{j,1}]_{v}-v[B_{i,k+1},B_{j}]_{v^{-1}}=0.$
		The identity in the proposition follows by applying $\TT_{\bome_j}^{-l}$ to this identity.
	\end{proof}

	\subsection{Relation \eqref{qsiDR2} when $c_{ij}=0=c_{\tau i,j}$}
	\label{subsec:verifydeg}
   We start with a general preparatory lemma. 
	\begin{lemma}
		\label{lem:comb}
		For $i\in\II$ and $n\ge 1$, we have
		\begin{align*}
			\Theta_{i,n}=
			\begin{cases}
				v^{-2} C \Theta_{i,n-2}+v^{2}\K_i^{-1} \big( [B_{ i,0},B_{i,n}]_{v^{-2}}+[B_{i,n-1},B_{ i,1}]_{v^{-2}} \big), & \text{for } n\ge 3,
				\\
				-v^{-1} C +v^{2}\K_i^{-1} \big( [B_{ i},B_{i,2}]_{v^{-2}}+[B_{i,1},B_{i,1}]_{v^{-2}} \big), & \text{for } n=2, \qquad \text{ if } \tau i=i,
				\\
				v^{2}\K_{i}^{-1} [B_{i}, B_{i,1}]_{v^{-2}}, & \text{for } n= 1,
			\end{cases}
		\end{align*}
		and $\Theta_{i,n}=[B_{i,n},B_{\tau i}]\K_{\tau i}^{-1}$
		if $c_{i,\tau i}=0$.
		In particular, the element $\Theta_{i,n}$  is a $\Q(v)[C^{\pm 1},\K_i^{\pm 1},\K_{\tau i}^{\pm1}]$-linear combination of $1$ and $[B_{i,k}, B_{\tau i,l+1}]_{v^{-c_{i,\tau i}}} +[B_{\tau i,l}, B_{i,k+1}]_{v^{-c_{i,\tau i}}}$, for $l, k \in \Z$.
	\end{lemma}
	
	\begin{proof}
 The formula for $\Theta_{i,n}$ when $c_{i,\tau i}=0$ is \eqref{qsiDR4} with $k=n$ and $l=0$. 
	The recursion formula in the lemma for $\Theta_{i,n}$ when $\tau i=i$ is a reformulation of \eqref{qsiDR5} with $k=n-1$ and $l=0$. The second statement follows by an induction on $n$ using the recursion formula. (A precise linear combination can be written down, but will not be needed.)
	\end{proof}
		
	Now we establish the relation \eqref{qsiDR2} under the assumption that $c_{ij}=0=c_{i,\tau j}$.
	\begin{proposition}
		\label{prop:iDR2}
		Assume $c_{ij}=0=c_{i,\tau j}$, for $i,j\in \II$. Then, for $m\geq 1$ and $r \in\Z$, we have
		\begin{align*}
			[\Theta_{i,m},B_{j,r}] &=0 =[H_{i,m},B_{j,r}].
		\end{align*}
	\end{proposition}

	\begin{proof}
		We shall prove the first equation only (as the two equations are clearly equivalent). By Lemma~\ref{lem:comb}, it suffices to check that $[B_{i,k}, B_{\tau i,l}]_{v^{-c_{i,\tau i}}}$ commutes with $B_{j,r}$ for all $k,l,r$. But this clearly follows from the identities $[ B_{i,k},B_{j,r}]=0=[B_{\tau i,l},B_{j,r}]$ by Proposition \ref{prop:qsiDR7}.
	\end{proof}

	\subsection{Relation \eqref{qsiDR2} when $(c_{ij}=-1, c_{\tau i,j}=0)$ or $(c_{ij}=0, c_{\tau i,j}=-1)$}

 The two cases are very similar, and we give the detailed proof of the relation  \eqref{qsiDR2} in case that $c_{ij}=-1$ and $c_{\tau i,j}=0$. In this case, we have $\tau i\neq i$, $\tau j\neq j$, and thus  $c_{i,\tau i}=0$. 
 \begin{proposition}
     Assume that $c_{ij}=-1$ and $c_{\tau i,j}=0$. Then $[H_{i,m}, B_{j,l}] = -\frac{[m]}{m} B_{j,m+l}$.
 \end{proposition}

 \begin{proof}
 By Lemma \ref{lem:HB}, the desired identity in the proposition is equivalent to the relation 
$\bB_j(w)  \bTH_i(z) = \frac{1 -v^{-1}zw^{-1}}{1 -vzw^{-1}} \bTH_i(z) \bB_j(w)$, which is then equivalent to the relation
	\begin{align}
		\label{eq:HBequiv1}
		[\Theta_{i,m},B_{j,l}]=v^{-1}[\Theta_{i,m-1},B_{j,l+1}]_{v^2}, \quad \text{ for }m\geq0 \text{ and }l\in\Z.
	\end{align}
By applying $\TT_{\bome_j}^{-l}$ and the invariance of $\Theta_{i,m}$ under $\TT_{\bome_j}$ in Proposition \ref{invTheta}, we reduce a proof of \eqref{eq:HBequiv1} to a proof of its special case at $l=0$:
\begin{align}
		\label{eq:HBequiv1b}
		[\Theta_{i,m},B_{j}]=v^{-1}[\Theta_{i,m-1},B_{j,1}]_{v^2}, \quad \text{ for }m\geq0.
	\end{align}

It remains to prove \eqref{eq:HBequiv1b}. We first deal with the case when $m=1$. Since $[B_{\tau i},B_j]=0$, by \eqref{def:thetac=0} and Lemma~ \ref{lem:fixB}, we have
	\begin{align*}
		[\Theta_{i,1},B_j]=	&o(i)^{-1}\big[[\TT_{\bome_i}^{-1}(B_i), B_{\tau i} ] \K_{\tau i}^{-1},B_j\big]\\
		&=o(i)^{-1} v\big[[\TT_{\bome_i}^{-1}(B_i), B_{\tau i} ] ,B_j\big]_{v^{-1}}\K_{\tau i}^{-1}
		\\
		&=o(i)^{-1} v\big[[\TT_{\bome_i}^{-1}(B_i), B_{j} ]_{v^{-1}} ,B_{\tau i}\big]\K_{\tau i}^{-1}
		\\
		&=-o(i)^{-1} \big[\TT_{\bome_i}^{-1}(X_{ij}) ,B_{\tau i}\big]\K_{\tau i}^{-1}.
	\end{align*}
	Using Lemmas \ref{lem:TXij} and \ref{lem:Tomeij}, we have
	\begin{align*}
		[\Theta_{i,1},B_j]
		&=-o(i)^{-1} \big[\TT_{\bome_j}^{-1}(X_{ji}) ,B_{\tau i}\big]\K_{\tau i}^{-1}
		\\
		&=-o(i)^{-1} \TT_{\bome_j}^{-1}\big([B_iB_j-vB_jB_i, ,B_{\tau i}]\K_{\tau i}^{-1}\big)
		\\
		&=-o(i)^{-1}\TT_{\bome_j}^{-1}\Big( \big[[B_i,B_{\tau i}],B_j\big]_v \K_{\tau i}^{-1}\Big)
		\\
		&= -o(i)^{-1}\TT_{\bome_j}^{-1}\big( [\frac{\K_{\tau i}-\K_i}{v-v^{-1}},B_j]\K_{\tau i}^{-1}\big)
		\\
		&=o(i)^{-1}\TT_{\bome_j}^{-1}(B_j)
		\\
		&=-B_{j,1}.
	\end{align*}
 This proves \eqref{eq:HBequiv1b} for $m=1$.
		
	For $m>0$, we have
	\begin{align*}
		[\Theta_{i,m},B_j]&= \big[[B_{i,m-1},B_{\tau i,-1}]\K_{\tau i}^{-1},B_j\big]C
		\\
		&=v\big[[B_{i,m-1},B_j]_{v^{-1}},B_{\tau i,-1}] \big]\K_{\tau i}^{-1}C
		\\
		&=o(i)^{1-m} v\big[[\TT_{\bome_i}^{1-m} (B_i),B_j]_{v^{-1}},B_{\tau i,-1}\big]\K_{\tau i}^{-1}C
		\\
		&=o(i)^{1-m} v\big[\TT_{\bome_i}^{1-m}([B_i,B_j]_{v^{-1}}),B_{\tau i,-1}\big]\K_{\tau i}^{-1}C
		\\
		&=-o(i)^{1-m} v\big[\TT_{\bome_i}^{2-m}T_{\bome_i}^{-1}(X_{ij}),B_{\tau i,-1}\big]\K_{\tau i}^{-1}C
		\\
		&=-o(i)^{1-m} v\big[\TT_{\bome_i}^{2-m}T_{\bome_j}^{-1}(X_{ji}),B_{\tau i,-1}\big]\K_{\tau i}^{-1}C
		\\
		&=-o(i) \TT_{\bome_i}^{2-m}T_{\bome_j}^{-1}\big[ [B_i,B_{\tau i ,1-m}],B_j\big]_v\K_{\tau i}^{-1}C.
	\end{align*}
	Using the relation \eqref{qsiDR4} which was already proved in Proposition \ref{prop:iDR31}, we have
	\begin{align*}
		[\Theta_{i,m},B_j]&=-o(i) \TT_{\bome_i}^{2-m}T_{\bome_j}^{-1}[\K_{\tau i}C^{1-m}\Theta_{i,m-1},B_j]_v \K_{\tau i}^{-1}C\\
		&=v^{-1}\Theta_{i,m-1}B_{j,1}-vB_{j,1}\Theta_{i,m-1}
		\\
		&=v^{-1} [\Theta_{i,m-1},B_{j,1}]_{v^2}.
	\end{align*}
	This proves \eqref{eq:HBequiv1b}.
 \end{proof}	

	\subsection{Relation \eqref{qsiDR2} when $c_{ij}=-1=c_{\tau i,j}$}

Assume that $c_{ij}=-1=c_{\tau i,j}$ in this subsection. 

In case $\tau i=i$, the proof of the relation \eqref{qsiDR2} is the same as in \cite[\S 4]{Z22} for split type, which only uses \eqref{qsiDR3} and Serre relation \eqref{eq:S2}. We will not repeat the argument here.

From now on we shall always assume that $\tau i\neq i$. 
Together with $c_{ij}=-1=c_{\tau i,j}$, this implies that $c_{i,\tau i}=0$ and $\tau j=j$ by inspection of Table \ref{tab:Satakediag}. By Lemma \ref{lem:HB}, the relation \eqref{qsiDR2} can be equivalently reformulated as follows.

\begin{proposition}
Assume that $\tau i\neq i$ and $c_{ij}=-1=c_{\tau i,j}$. Then the following relation holds,  for $m \ge 0 \text{ and  } r \in \Z$:
	\begin{align}
		[ & \Theta_{i,m},B_{j,r}]+[\Theta_{i,m-2},B_{j,r}]C
		\notag \\
		&= v^{-1}[\Theta_{i,m-1},B_{j,r+1}]_{v^2}+v[\Theta_{i,m-1},B_{j,r-1}]_{v^{-2}}C.
		\label{iDR2-reform}
	\end{align}
\end{proposition}

\begin{proof}
By the invariance of $\Theta_{i,m}$ under $\TT_{\bome_j}$ in Proposition \ref{invTheta} and applying $\TT_{\bome_j}^{r}$ to \eqref{iDR2-reform}, we have reduced its proof to a proof of the special case \eqref{iDR2-reform}$|_{r=0}$. 

 We first prove the identity \eqref{iDR2-reform}$|_{r=0}$ for $m=1$. 
Note also that $o(i)=o(\tau i)$ and then $\btau \Theta_{i,1} =\Theta_{\tau i,1}$. We have
	\begin{align*}
		[\Theta_{i,1},B_j]&=o(i)^{-1}\big[ [\TT_{\bome_i}^{-1}(B_i),B_{\tau i}]\K_{\tau i}^{-1},B_j \big]
		\\
		&=o(i)^{-1}\big[ [\TT_{\bome_i}^{-1}(B_i),B_{\tau i}],B_j \big]\K_{\tau i}^{-1}
		\\
		&=o(i)^{-1}v\big[ [\TT_{\bome_i}^{-1}(B_i),B_{j}]_{v^{-1}},B_{\tau i} \big]_{v^{-1}}\K_{\tau i}^{-1}
		-o(i)^{-1}v^{-1}\big[ [B_{\tau i},B_{j}]_{v},\TT_{\bome_i}^{-1}(B_{i}) \big]_{v}\K_{\tau i}^{-1}
		\\
		&=-o(i)^{-1}[\TT_{\bome_i}^{-1}(X_{ij}),B_{\tau i}]_{v^{-1}}\K_{\tau i}^{-1} -o(i)^{-1}v^{-1}\TT_{\bome_i}^{-1} \big([\TT_{\bome_i}(X_{j,\tau i}),B_i]_v\big)\K_{\tau i}^{-1}
		\\
		&=-o(i)^{-1}[\TT_{\bome_j}^{-1}(X_{ji}),B_{\tau i}]_{v^{-1}}\K_{\tau i}^{-1} -o(i)^{-1}v^{-1}\TT_{\bome_i}^{-1} \big([\TT_{\bome_j}(X_{\tau i,j}),B_i]_v\big)\K_{\tau i}^{-1}
		\\
		&=-o(i)^{-1}\TT_{\bome_j}^{-1}\Big( \big[[B_i,B_j]_v,B_{\tau i}\big]_{v^{-1}}\Big)\K_{\tau i}^{-1} -o(i)^{-1}v^{-1}\TT_{\bome_i}^{-1}\TT_{\bome_j} \Big(\big[[B_j,B_{\tau i}]_v,B_i\big]_v\Big)\K_{\tau i}^{-1}
		\\
		&=\big[[B_i,B_{j,1}]_v,B_{\tau i}\big]_{v^{-1}}\K_{\tau i}^{-1}
		+v^{-1}\big[[B_{j,-1},B_{\tau i,1}]_v,B_{i,1}\big]_v\K_{\tau i}^{-1}
		\\
		&=\big( B_iB_{j,1}B_{\tau i}-vB_{j,1}B_iB_{\tau i}-v^{-1}B_{\tau i}B_iB_{j,1}+B_{\tau i}B_{j,1}B_i \big) \K_{\tau i}^{-1}
		\\
		&+\big( v^{-1} B_{j,-1}B_{\tau i,1}B_{i,1}-B_{\tau i,1}B_{j,-1}B_{i,1}-B_{i,1}B_{j,-1}B_{\tau i,1} +vB_{i,1}B_{\tau i,1}B_{j,-1} \big)\K_{\tau i}^{-1}.
	\end{align*}
	
By Theorem~\ref{thm:Ti}, we have 
\begin{align*}
    \TT_i(B_j)=\TT_{\tau i}(B_j) &=\big[[B_j,B_{\tau i}]_v,B_{i}\big]_v-v B_j \K_{\tau i}
    \\
    \TT_i^{-1}(B_j)=\sigma_\tau \TT_i (B_j) &=\big[B_{\tau i},[B_i, B_j]_v\big]_v-v B_j \K_{\tau i}.
\end{align*}
Since $\bs_i t_{\bome_j}=t_{\bome_j}\bs_i$, we have $\TT_i \TT_{\bome_j}=\TT_{\bome_j} \TT_i$. Then, by applying $o(j)\TT_{\bome_j}$ to the two identities above, we have
	\begin{align}
		\label{eq:pfiDR2'}
		\begin{split}
			&\TT_i(B_{j,-1})=\big[[B_{j,-1},B_{\tau i}]_v,B_{i}\big]_v-v B_{j,-1} \K_{\tau i},
   \\
			&\TT_i^{-1}(B_{j,1})=\big[B_{\tau i},[B_i, B_{j,1}]_v\big]_v-v B_{j,1} \K_{\tau i}.
		\end{split}
	\end{align}
	Using \eqref{eq:pfiDR2'}, we rewrite the above formula of $[\Theta_{i,1},B_j]$ as follows:
	\begin{align*}
		[\Theta_{i,1},B_j]
		&=\big[[B_i,B_{j,1}]_v,B_{\tau i}\big]_{v^{-1}}\K_{\tau i}^{-1}
		+v^{-1}\big[[B_{j,-1},B_{\tau i,1}]_v,B_{i,1}\big]_v\K_{\tau i}^{-1}
		\\
		&=-v^{-1}\big[B_{\tau i},[B_i,B_{j,1}]_v\big]_{v }\K_{\tau i}^{-1}
		+v^{-1}\big[[B_{j,-1},B_{\tau i,1}]_v,B_{i,1}\big]_v\K_{\tau i}^{-1}
		\\
		&=-v^{-1} \big(\TT_i^{-1}(B_{j,1}) +v B_{j,1} \K_{\tau i}\big)\K_{\tau i}^{-1} + v^{-1} \TT_{\bome_i}^{-1}\big(\TT_i(B_{j,-1})+v B_{j,-1} \K_{\tau i}\big)\K_{\tau i}^{-1}
		\\
		&=-B_{j,1}+B_{j,-1} C -v^{-1}\TT_i^{-1}(B_{j,1})\K_{\tau i}^{-1} +v^{-1}\TT_{\bome_i}^{-1}\TT_i(B_{j,-1})\K_{\tau i}^{-1}
	\end{align*}
	By Lemma~\ref{lem:Tomeij}, we have $\TT_{\bome_i}^{-1}\TT_i(B_{j,-1})=\TT_i^{-1}\TT_{\bome_i}\TT_{\bome_j}^{-2}(B_{j,-1})=\TT_i^{-1}(B_{j,1})$. Then we simplify the previous identity to be 	$[\Theta_{i,1},B_j]=-B_{j,1}+B_{j,-1} C$, 
which is the desired identity \eqref{iDR2-reform}$|_{r=0}$ for $m=1.$
	
Now let $m\geq2$. By \eqref{def:thetac=0} we have
$\Theta_{i,m}=[B_{i,m},B_{\tau i}]\K_{\tau i}^{-1}.$ Then we have
	\begin{align*}
		[\Theta_{i,m},B_j]\K_{\tau i} &= \big[[B_{i,m},B_{\tau i}] ,B_j\big]
		\\
		&= \big[B_{i,m},[B_{\tau i},B_j]_v\big]_{v^{-1}} -\big[ B_{\tau i},[B_{i,m},B_j]_{v^{-1}} \big]_v
		\\
		&=v\big[ B_{i,m},[B_{\tau i,1},B_{j,-1}]_{v^{-1}} \big]_{v^{-1}} -v^{-1} \big[ B_{\tau i}, [B_{i,m-1},B_{j,1}] \big],
	\end{align*}
	where the last equality follows from \eqref{qsiDR3} which has been proved. Therefore, we have
	\begin{align}
		[\Theta_{i,m},B_j]\K_{\tau i} &=
		-\big[B_{\tau i,1},[B_{j,-1},B_{i,m}]_v\big]_{v^{-1}}+v^{-1}\big[B_{j,-1},[B_{\tau i,1},B_{i,m}]\big]_{v^2}
  \notag
		\\\notag
		&\quad +\big[B_{i,m-1},[B_{j,1},B_{\tau i}]_{v^{-1}}\big]_{v}-v\big[B_{j,1},[B_{i,m-1},B_{\tau i}]\big]_{v^{-2}}
		\\\notag
		&=	-v\big[B_{\tau i,1},[B_{j},B_{i,m-1}]_{v^{-1}}\big]_{v^{-1}}-v^{-1}\big[B_{j,-1},\Theta_{i,m-1}\big]_{v^2}C\K_{\tau i}
		\\\notag
		&\quad +v^{-1}\big[B_{i,m-1},[B_{j},B_{\tau i,1}]_{v}\big]_{v}-v\big[B_{j,1},\Theta_{i,m-1}\big]_{v^{-2}}\K_{\tau i}
		\\
		&=\big[B_j,[B_{i,m-1},B_{\tau i,1}]\big]+v^{-1}[\Theta_{i,m-1},B_j]_{v^2}\K_{\tau i}+v[\Theta_{i,m-1},B_{j,-1}]_{v^{-2}}\K_{\tau i}, 		\label{eq:proofThetaB}
	\end{align}
	where the second equality follows from \eqref{qsiDR3}.
	If $m=2$, then we have \begin{align*}
		\big[B_j,[B_{i,m-1},B_{\tau i,1}]\big]&=\big[B_j,\frac{\K_{\tau i}C-\K_{i}C}{v-v^{-1}}\big]=0.
	\end{align*}
	If $m\geq2$, we have
	\begin{align*}
		\big[B_j,[B_{i,m-1},B_{\tau i,1}]\big]&=[B_j,\Theta_{i,m-2}]C\K_{\tau i}.
	\end{align*}
	Plugging this into \eqref{eq:proofThetaB}, we have proved the desired identity \eqref{iDR2-reform}$|_{r=0}$ and hence the proposition. 
\end{proof}

	\subsection{Relation \eqref{qsiDR1} for $i\neq j$}
	
	In this subsection, we shall derive the identity \eqref{qsiDR1} (i.e., $[H_{i,m},H_{j,n}]=0$, for $i\neq j \in \II$) from the relations \eqref{qsiDR2}--\eqref{qsiDR5} which have been established above. We start with some preparatory lemma. 
	
	\begin{lemma}
		\label{lem:comm}
		For $i,j \in \II$, $l, k \in \Z$ and  $m\ge 1$, we have
		\[
		\Big[ H_{i,m}, [B_{j,k},B_{\tau j,l+1}]_{v^{-c_{j,\tau j}}}+[B_{\tau j,l},B_{j,k+1}]_{v^{-c_{j,\tau j}}} \Big] =0.
		\]
	\end{lemma}
	
	\begin{proof}
		If $j=\tau j$, then it follows from the same proof of \cite[Lemma 4.7]{LW21b} by using \eqref{qsiDR5}.
		
		For $j\neq \tau j$, using \eqref{qsiDR2} we have
		\begin{align*}
		\Big[ H_{i,m},  &  [B_{j,k},B_{\tau j,l+1}]_{v^{-c_{j,\tau j}}}+[B_{\tau j,l},B_{j,k+1}]_{v^{-c_{j,\tau j}}} \Big]
			\\\notag
			=&\, [H_{i,m}, B_{j,k}] B_{\tau j,l+1} +v^{-c_{i,\tau j}} B_{\tau j,l+1} [B_{j,k}, H_{i,m}]
			\\\notag
			& + [H_{i,m}, B_{\tau j,l}] B_{j,k+1} +v^{-c_{i,\tau j}} B_{j,k+1} [B_{\tau j,l}, H_{i,m}]
			\\\notag
			& +B_{\tau j,l} [H_{i,m}, B_{j,k+1}]  +v^{-c_{i,\tau j}}  [B_{j,k+1}, H_{i,m}] B_{\tau j,l}
			\\\notag
			& + B_{j,k} [H_{i,m}, B_{\tau j,l+1}]  +v^{-c_{i,\tau j}}  [B_{\tau j,l+1}, H_{i,m}] B_{j,k}\\\notag
		=&\, \Big(\frac{[mc_{ij}]}{m} B_{j, k+m} -\frac{[mc_{i,\tau j}]}{m}B_{j,k-m} C^m\Big) B_{\tau j, l+1}
			\\
			&-v^{-c_{i,\tau j}} B_{\tau j,l+1} \Big(\frac{[mc_{ij}]}{m} B_{j,k+m} -\frac{[mc_{i,\tau j}]}{m}B_{j,k-m} C^m\Big)
			\\\notag
			& + \Big(\frac{[mc_{i,\tau j}]}{m}B_{\tau j, l+m} -\frac{[mc_{ij}]}{m}B_{\tau j,l-m} C^m\Big) B_{j, k+1}
			\\
			&-v^{-c_{i,\tau j}} B_{j,k+1} \Big(\frac{[mc_{i,\tau j}]}{m} B_{\tau j,l+m} -\frac{[mc_{ij}]}{m}B_{\tau j,l-m} C^m\Big)
			\\\notag
			& + B_{\tau j, l}\Big(\frac{[mc_{ij}]}{m}B_{j, k+m+1} -\frac{[mc_{i,\tau j}]}{m}B_{j,k-m+1} C^m\Big)
			\\
			&-v^{-c_{i,\tau j}} \Big(\frac{[mc_{ij}]}{m} B_{j,k+m+1} -\frac{[mc_{i,\tau j}]}{m}B_{j,k-m+1} C^m\Big) B_{\tau j, l}
			\\\notag
			& + B_{j, k}\Big(\frac{[mc_{i,\tau j}]}{m} B_{\tau j, l+m+1} -\frac{[mc_{ij}]}{m}B_{\tau j,l-m+1} C^m\Big)
			\\
			&-v^{-c_{i,\tau j}} \Big(\frac{[mc_{i,\tau j}]}{m}B_{\tau j,l+m+1} -\frac{[mc_{ij}]}{m}B_{\tau j,l-m+1} C^m\Big) B_{j, k}.
		\end{align*}
		The above equality can be further rewritten as
		\begin{align}
			\Big[ H_{i,m}, &[B_{j,k},B_{\tau j,l+1}]_{v^{-c_{j,\tau j}}}+[B_{\tau j,l},B_{j,k+1}]_{v^{-c_{j,\tau j}}} \Big]
            \notag \\
            = &\, \frac{[mc_{ij}]}{m}\Big([B_{j,k+m}, B_{\tau j,l+1}]_{v^{-c_{i,\tau j}}} + [B_{\tau j,l}, B_{j,k+m+1}]_{v^{-c_{i,\tau j}}}\Big)
			\notag \\
			&-\frac{[mc_{ij}]}{m} \Big([B_{j,k}, B_{\tau j,l-m+1}]_{v^{-c_{i,\tau j}}}  C^m +[B_{\tau j,l-m}, B_{ j,k+1}]_{v^{-c_{i,\tau j}}}  C^m\Big)
			\notag \\
			&+\frac{[mc_{i,\tau j}]}{m}\Big([B_{j,k}, B_{\tau j,l+m+1}]_{v^{-c_{i,\tau j}}}+ [B_{\tau j,l+m}, B_{j,k+1}]_{v^{-c_{i,\tau j}}} \Big) 
			\notag \\
			&-\frac{[mc_{i,\tau j}]}{m}\Big( [B_{j,k-m}, B_{\tau j,l+1}]_{v^{-c_{i,\tau j}}}  C^m + [B_{\tau j,l}, B_{j,k-m+1}]_{v^{-c_{i,\tau j}}}  C^m\Big).
   \label{HBB4}
		\end{align}

  Note that $c_{j,\tau j}=0$ since $j\neq \tau j$. Applying the identity \eqref{qsiDR4} repeatedly gives us
		\begin{align*}
			[B_{j,k+m}, B_{\tau j,l+1}]_{v^{-c_{i,\tau j}}}=[B_{j,k}, B_{\tau j,l-m+1}]_{v^{-c_{i,\tau j}}}  C^m,
			\\
			[B_{\tau j,l}, B_{j,k+m+1}]_{v^{-c_{i,\tau j}}}=[B_{\tau j,l-m}, B_{ j,k+1}]_{v^{-c_{i,\tau j}}}  C^m.
		\end{align*}
  This implies that the first two summands involving $\frac{[mc_{ij}]}{m}$ on the RHS of \eqref{HBB4} cancel out. Similarly the last two summands involving $\frac{[mc_{i,\tau j}]}{m}$ on the RHS of \eqref{HBB4} also cancel out. This proves that the LHS of \eqref{HBB4} $ =0$ as desired. 
	\end{proof}
We can now establish the relation \eqref{qsiDR1} for $i\neq j  \in \II$.
	\begin{proposition}
		\label{prop:iDR1}
		The identity $[H_{i,m},H_{j,n}]=0$, for $i\neq j  \in \II$, holds in $\tUi$.
	\end{proposition}
	
	\begin{proof}
It follows by Lemma~\ref{lem:comb} and Lemma~\ref{lem:comm} that $[ H_{i,m}, \Theta_{j,a}] =0$, for all $m, a\ge 1.$ Since $H_{j,m}$ is a linear combination of monomials in $\Theta_{j,b}$'s by \eqref{exp}, we conclude that $[ H_{i,m}, H_{j,n}] =0$.
	\end{proof}

	\subsection{Relations \eqref{qsiDR9}--\eqref{qsiDR8}}
	
	We shall fix $i,j \in \II$ such that $c_{ij}=-1$ throughout this subsection.
	
	\begin{lemma}
		\label{lem:SSSa}
		For $k_1, k_2, l \in \Z$, we have
		\begin{align*}
			& \SS(k_1,k_2+1 |l) + \SS(k_1+1,k_2|l) -[2] \SS(k_1+1,k_2+1 |l-1)
			\\
			&=\delta_{i,\tau i} \Big( -[B_{jl}, \Theta_{i, k_2-k_1+1}]_{v^{-2}} C^{k_1} +v^{-2} [B_{jl}, \Theta_{i, k_2-k_1-1}]_{v^{-2}} C^{k_1+1} \Big) \K_i + \{k_1 \leftrightarrow k_2\}.
		\end{align*}
	\end{lemma}

	\begin{proof}
		The proof is the same as for \cite[Lemma 4.9]{LW21b}, and hence omitted here. It uses only the relations \eqref{qsiDR5}--\eqref{qsiDR3}, which have been established above.
	\end{proof}

	\begin{lemma}  \label{lem:SSS}
		For $k_1, k_2, l \in \Z$, we have
		\begin{align*}
			& \SS(k_1,k_2+1 |l) + \SS(k_1+1,k_2|l) -[2] \SS(k_1,k_2|l+1)
			\\
			&= \delta_{i,\tau i}\Big( -[\Theta_{i, k_2-k_1+1}, B_{jl}]_{v^{-2}} C^{k_1} +v^{-2} [\Theta_{i, k_2-k_1-1}, B_{jl}]_{v^{-2}} C^{k_1+1} \Big) \K_i + \{k_1 \leftrightarrow k_2\}.
		\end{align*}
	\end{lemma}
	
	\begin{proof}
		The proof is the same as for \cite[Lemma 4.13]{LW21b}, and hence omitted here.	It uses only the relations \eqref{qsiDR5}--\eqref{qsiDR3}.
	\end{proof}
	
	\begin{proposition}
		Relations \eqref{qsiDR9}--\eqref{qsiDR8} hold in $\tUi$.
	\end{proposition}
	
	\begin{proof}
		Follows by the same argument of \cite[\S5.1]{Z22} with the help of Lemmas \ref{lem:SSSa}--\ref{lem:SSS}.
	\end{proof}	 
 
Summarizing the results in this section, we have verified all the relations \eqref{qsiDR0}--\eqref{qsiDR8} of $\tUiD$ for the images of the generators under $\Phi: \tUiD \rightarrow\tUi$. This completes the proof of Theorem \ref{thm:Dr}.

\appendix

\section{Proofs of Lemma \ref{lem:weight} and Lemma \ref{lem:qsweight}}  
 
In this appendix, we provide the detailed proofs for Lemma \ref{lem:weight} and Lemma \ref{lem:qsweight}.

\subsection{Proof of Lemma~\ref{lem:weight} }
\label{app:A1}

 \begin{proof}
 For convenience we recall Lemma~\ref{lem:weight} states that, for $\tau i =i\in \II$, we have
\begin{align}
		\label{eq:embedding}
		&\Big[\big[[ \TT_{\omega'_i}^{-1}(B_i) ,B_i]_{v^2} , B_i\big],B_i\Big]_{v^{-2}}-v[2]^2[\TT_{\omega'_i}^{-1}(B_i),B_i] \K_i=0.
	\end{align}
Let us verify the identity \eqref{eq:embedding}.
  
We first recall several ingredients from \cite{WZ23}, which will be used in the proof.
Let $\tcT_i$ and $\tcT_w$ be the (rescaled) Lusztig braid group symmetries on $\tU$ introduced in \cite[\S 4.1]{WZ23}. Let $\tfX$ be the quasi $K$-matrix for the quantum symmetric pair $(\tU,\tUi)$ and $\tfX_i$ be rank one quasi $K$-matrices, for $i\in \I$.
    \cite[Theorem 3.6]{WZ23} is applicable to the $\TT_i$ in Theorem~\ref{thm:Ti} in the split affine type, and it gives us
    \begin{align}  \label{intertwiner}
    \TT_i^{-1}(x) \tfX_i=\tfX_i\tcT_{i}^{-1}(x),
    \end{align} 
    for any $x\in \tUi$.

 By assumption $\tau i=i$, we have $\bs_i=s_i$, $t_{\bome_i} =t_{\omega_i}$ and $t_{\bome_i'} =t_{\omega_i'}$. Fix a reduced expression $t_{\omega'_i} =\tau_i s_{i_1}s_{i_2} \cdots s_{i_k}$ where $\tau_i$ is a Dynkin diagram automorphism. We define
	\begin{align*}
		\tfX_{\omega'_i}=\tau_i \tfX_{i_1} \tcT'_{i,-1}(\tfX_{i_2})\cdots \tcT'_{s_{i_1}s_{i_2}\cdots s_{i_{k-1}},-1}(\tfX_{i_k}).
	\end{align*}
	($\tfX_{\omega'_i}$ here differs from the element in \cite[(8.1)]{WZ23} by the anti-involution $\sigma$. By \cite[Theorem 8.1]{WZ23} $\tfX_{\omega'_i}$ is independent of the choice of a reduced expression of $t_{\omega'_i}$.) Then a repeated use of \eqref{intertwiner} gives us
	\begin{align}\label{eq:embedding3}
		\TT_{\omega'_i}^{-1}(B_i) \tfX_{\omega'_i}
		=\tfX_{\omega'_i}\tcT'_{\omega'_i,-1}(B_i).
	\end{align}

    To verify  \eqref{eq:embedding}, we shall apply a strategy similar to \cite[Proof of Proposition 7.2]{WZ23} (where the notation $\tcT'_{w,-1} =\tcT_w^{-1}$ is used). Let us write the LHS of \eqref{eq:embedding} in terms of a monomial basis of $\tUi$ as
	\begin{align}
			\label{eq:mbasis}
		&\Big[\big[[ \TT_{\omega'_i}^{-1}(B_i) ,B_i]_{v^2} , B_i\big],B_i\Big]_{v^{-2}}-v[2]^2[\TT_{\omega'_i}^{-1}(B_i),B_i] \K_i=\sum_{J} a_J B_J,
	\end{align}
	where $B_J=B_{j_1} B_{j_2} \cdots B_{j_k}$ for some multi-index $J=(j_1,j_2,\ldots,j_k)$ with $j_a\in\I$ and $a_J \in \U^{\imath 0}$.  The lowest weight homogeneous component of $B_J$ (viewed as an element in $\tU$) has weight $-\mathrm{wt}(J)=-\alpha_{j_1} -\ldots-\alpha_{j_k}$. Denote by $Q_i:=\{-\mathrm{wt}(J)\mid a_J\neq 0\}$ the set of negative weights whose corresponding homogeneous components of the LHS of \eqref{eq:mbasis} are nonzero. We shall prove that $Q_i=\emptyset.$
		
Recall from \eqref{eq:embedding3} that
$\TT_{\omega'_i}^{-1}(B_i)
=\tfX_{\omega'_i}\tcT'_{\omega'_i,-1}(B_i)\tfX_{\omega'_i}^{-1}$,
and by \cite[Lemma~ 3.6]{DK19}, we have that $\tfX_{\omega'_i}\in \bigoplus_{ \alpha\in \N\Delta(\omega_i')} \tU^+_{\alpha+\tau \alpha}$. Then any negative weight appearing in the weights of homogeneous components of $\TT_{\omega'_i}^{-1}(B_i)$ must have the form $\alpha_i-\delta+\gamma+\tau \gamma$, for some $\gamma\in \N\Delta(\omega_i')$. Moreover, since $\tau i=i$, the set $\Delta(\omega_i')$ is fixed by $\tau$. Then, by Lemma~\ref{C:weight}, the multiplicity of $\alpha_0$ in $\gamma+\tau \gamma$ for nonzero $\gamma\in\N\Delta(\omega_i')$ is at least $2$, and recall that the multiplicity of $\alpha_0$ in $\delta$ is 1. Hence, by counting the multiplicity of $\alpha_0$, the only negative weight appearing in homogeneous components of $\TT_{\omega'_i}^{-1}(B_i)$ is $\alpha_i-\delta$. It follows by definition that $Q_i \subseteq \{-2\alpha_i-\delta,-\delta,2\alpha_i-\delta\}$. 

A direct computation shows that the homogeneous components of \eqref{eq:embedding} of the weights $-2\alpha_i-\delta, -\delta, 2\alpha_i-\delta$ are the LHS of \eqref{eq:iBeck1}--\eqref{eq:iBeck3} below, which we claim are all equal to 0:
	\begin{align}\label{eq:iBeck1}
		\Big[\big[[ \tcT'_{\omega'_i,-1}(F_i) ,F_i]_{v^2} , F_i\big],F_i\Big]_{v^{-2}}&=0,
			\\\notag
			\Big[\big[[ \tcT'_{\omega'_i,-1}(F_i) ,F_i]_{v^2} , F_i\big],E_i K_i'\Big]_{v^{-2}}
			&+\Big[\big[[ \tcT'_{\omega'_i,-1}(F_i) ,F_i]_{v^2} , E_i K_i'\big],F_i\Big]_{v^{-2}}\\
			+\Big[\big[[ \tcT'_{\omega'_i,-1}(F_i) ,E_i K_i']_{v^2} ,F_i\big],F_i\Big]_{v^{-2}} &-v[2]^2[ \tcT'_{\omega'_i,-1}(F_i) ,F_i] K_i K_i' =0,
			\label{eq:iBeck2}
			\\\notag
			\Big[\big[[ \tcT'_{\omega'_i,-1}(F_i) ,F_i]_{v^2} ,E_i K_i'\big],E_i K_i'\Big]_{v^{-2}}
			&+\Big[\big[[ \tcT'_{\omega'_i,-1}(F_i) ,E_i K_i']_{v^2} , E_i K_i'\big],F_i\Big]_{v^{-2}} \\
			+\Big[\big[[ \tcT'_{\omega'_i,-1}(F_i) ,E_i K_i']_{v^2} ,F_i\big],E_i K_i'\Big]_{v^{-2}} &-v[2]^2[ \tcT'_{\omega'_i,-1}(F_i) ,E_i K_i'] K_i K_i' =0.
			\label{eq:iBeck3}
		\end{align}
Indeed, the identity \eqref{eq:iBeck1} follows by \cite[Proposition 3.7]{Be94}. The identities \eqref{eq:iBeck2}--\eqref{eq:iBeck3} follow from the identity $\tcT'_{\omega'_i,-1}(F_i) \; E_i  = E_i  \; \tcT'_{\omega'_i,-1}(F_i)$ (see \cite[Lemma~ 3.4]{Be94}) and $[E_i,F_i]=\frac{K_i-K_i'}{v-v^{-1}}$. 

This shows that $Q_i=\emptyset$ and hence the RHS of \eqref{eq:mbasis} $=0$. This proves the desired identity \eqref{eq:embedding}. 
\end{proof}

\subsection{Proof of Lemma~\ref{lem:qsweight} }
\label{app:A2}

\begin{proof} 
It amounts to proving the identities \eqref{eq:qsembedding1}, \eqref{eq:qsembedding1a} and \eqref{eq:qsembedding2}. 
The identity  $[B_i,B_{\tau i}] =\frac{\K_{\tau i}-\K_i}{v-v^{-1}}$ in \eqref{eq:qsembedding2} follows by a direct computation. 

For convenience we recall the identities \eqref{eq:qsembedding1}--\eqref{eq:qsembedding1a} here as
\begin{align}
   \label{eq:qsembed1}
		\Big[\big[[ \TT_{\bome'_i}^{-1}(B_i) ,B_i]_{v^2}, B_i\big],B_i\Big]_{v^{-2}} &=0,\\
		[\TT_{\bome'_i}^{-1}(B_i),B_{\tau i}] &=0,
  \label{eq:qsembed1a}
	\end{align}
 for $i\in\II$ with $c_{i,\tau i}=0$. Let us verify \eqref{eq:qsembed1}--\eqref{eq:qsembed1a}; the overall idea is similar to the proof of the identity \eqref{eq:embedding}. 
 
Recall the braid group symmetry $\TT_{i}^{-1}$ (denoted by $\TT'_{i,-1}$ in \cite[Theorem 3.6]{WZ23}) satisfies that 
\begin{align} \label{intertwiner2}
    \TT'_{i,-1}(x) \tfX_i=\tfX_i\tcT'_{\bs_i,-1}(x), 
\end{align}
for any $x\in \tUi$. For a 
reduced expression $t_{\bome_i'}=\tau_i \bs_{i_1}\bs_{i_2} \cdots \bs_{i_k}$ with $\tau_i$ a diagram automorphism, we define
	\begin{align*}
		\tfX_{\bome'_i}=\tau_i \tfX_{\bs_{i_1}} \tcT'_{i,-1}(\tfX_{\bs_{i_2}})\cdots \tcT'_{\bs_{i_1}\bs_{i_2}\cdots \bs_{i_{k-1}},-1}(\tfX_{i_k}).
	\end{align*}
	By \cite[Theorem 8.1]{WZ23}, $\tfX_{\bome'_i}$ is independent of the choice of a reduced expression of $ t_{\bome_i'} $. Then a repeated application of \eqref{intertwiner2} gives us
	\begin{align}  \label{intetwiner3}
		\TT_{\bome'_i}^{-1}(B_i) \tfX_{\bome'_i}
		=\tfX_{\bome'_i}\tcT'_{\bome'_i,-1}(B_i).
	\end{align}

We prove \eqref{eq:qsembed1}. Write the LHS of \eqref{eq:qsembed1} in terms of a monomial basis of $\tUi$ as 
	\begin{align}
			\label{eq:qsmbasis}
		&\Big[\big[[ \TT_{\bome'_i}^{-1}(B_i) ,B_i]_{v^2} , B_i\big],B_i\Big]_{v^{-2}}=\sum_{J} a_J B_J,
	\end{align}
where $B_J=B_{j_1} B_{j_2} \cdots B_{j_k}$ for a multi-index $J=(j_1,j_2,\ldots,j_k)$ with $j_a\in\I$ and coefficients $a_J$. The lowest weight homogeneous component of $B_J$ (viewed as an element in $\tU$) has weight $-\mathrm{wt}(J)=-\alpha_{j_1} -\ldots-\alpha_{j_k}$. Denote by $Q_i:=\{-\mathrm{wt}(J)\mid a_J\neq 0\}$ the set of negative weights (in $-\N\I$) whose corresponding homogeneous components of the LHS of \eqref{eq:qsmbasis} are nonzero. We will show that $Q_i =\emptyset.$
		 
Note from \eqref{intetwiner3} that
$\TT_{\bome'_i}^{-1}(B_i)
=\tfX_{\bome'_i}\tcT'_{\bome'_i,-1}(B_i)\tfX_{\bome'_i}^{-1}$,
and by \cite[Lemma~ 3.6]{DK19}, we have $\tfX_{\bome'_i}\in \bigoplus_{ \balpha\in \N\Delta(\bome_i')} \tU^+_{2\balpha }$. Then any negative weight component for $\TT_{\bome'_i}^{-1}(B_i)$ has weight $\alpha_i-\delta+2\ov{\gamma}$, for some $\ov{\gamma}\in \N\Delta(\bome_i')$. Moreover, by Lemma~\ref{C:qsweight}, the multiplicity of $\alpha_0=\balpha_0$ in $2\ov{\gamma}$ for $\ov{\gamma}\in\N\Delta(\bome_i')\setminus\{0\}$ is at least $2$, and recall that the multiplicity of $\alpha_0$ in $\delta$ is 1. Hence, by counting the multiplicity of $\alpha_0$, the only negative weight component for $\TT_{\bome'_i}^{-1}(B_i)$ has weight $\alpha_i-\delta$.
It follows by definition that $Q_i\subset \{-2\alpha_i-\delta, -\alpha_i-\delta+\alpha_{\tau i}, 2\alpha_{\tau i}-\delta\}$. 
	
A direct computation shows that the homogeneous components of \eqref{eq:qsembed1} of the weights $-2\alpha_i-\delta, -\alpha_i-\delta+\alpha_{\tau i}, 2\alpha_{\tau i}-\delta$ are the LHS of \eqref{eq:qsiBeck1}--\eqref{eq:qsiBeck3} below, which we claim are all equal to 0:
		\begin{align}
			\label{eq:qsiBeck1}
			&\Big[\big[[ \tcT'_{\bome'_i,-1}(F_i) ,F_i]_{v^2} , F_i\big],F_i\Big]_{v^{-2}}=0,\\
			\notag
			&\Big[\big[[ \tcT'_{\bome'_i,-1}(F_i) ,E_{\tau i}K_i']_{v^2} , F_i\big],F_i\Big]_{v^{-2}}+\Big[\big[[ \tcT'_{\bome'_i,-1}(F_i) ,F_i]_{v^2} , E_{\tau i}K_i'\big],F_i\Big]_{v^{-2}}\\
			&+\Big[\big[[ \tcT'_{\bome'_i,-1}(F_i) ,F_i]_{v^2} , F_i\big],E_{\tau i}K_i'\Big]_{v^{-2}}=0,
   \label{eq:qsiBeck2}
			\\\notag
			&\Big[\big[[ \tcT'_{\bome'_i,-1}(F_i) ,E_{\tau i}K_i']_{v^2} , E_{\tau i}K_i'\big],F_i\Big]_{v^{-2}}+\Big[\big[[ \tcT'_{\bome'_i,-1}(F_i) ,F_i]_{v^2} , E_{\tau i}K_i'\big],E_{\tau i}K_i'\Big]_{v^{-2}}\\
			&+\Big[\big[[ \tcT'_{\bome'_i,-1}(F_i) ,E_{\tau i}K_i']_{v^2} , F_i\big],E_{\tau i}K_i'\Big]_{v^{-2}}=0.
   \label{eq:qsiBeck3}
		\end{align}	
Indeed, the identity \eqref{eq:qsiBeck1} follows by \cite[Proposition 3.7]{Be94}. The identities \eqref{eq:qsiBeck2}--\eqref{eq:qsiBeck3} can be derived from $[ E_{\tau i},F_i]=0$ and  $\tT'_{\bome'_i,-1}(F_i) E_{\tau i}= E_{\tau i} \tT'_{\bome'_i,-1}(F_i)$; see \cite[Lemma 3.5]{Be94}.

This proves that $Q_i =\emptyset$, i.e., the RHS of \eqref{eq:qsmbasis} $=0$, and hence the identity \eqref{eq:qsembed1}.
		
Just as the identity \eqref{eq:qsembed1} was reduced to the identities \eqref{eq:qsiBeck1}--\eqref{eq:qsiBeck3} above, the identity \eqref{eq:qsembed1a} is proved by an entirely similar argument (we skip the details) by reducing to the following identities \eqref{Be1}--\eqref{Be2}:
\begin{align}
  [\tcT'_{\bome'_i,-1}(F_i),E_{i}K_{\tau i}'] &=0,
   \label{Be1}\\
		[\tcT'_{\bome'_i,-1}(F_i),F_{\tau i}] &=0.
   \label{Be2}
\end{align}
The identity \eqref{Be1} follows from \cite[Lemma 3.4]{Be94}, while \eqref{Be2} follows by a direct computation
$[\tcT'_{\bome'_i,-1}(F_i),F_{\tau i}]
=\tcT'_{\bome'_i,-1}\big([F_i,F_{\tau i}]\big)=0$, thanks to $\tcT'_{\bome'_i,-1}(F_{\tau i})=F_{\tau i}$.
	
This completes the proof of Lemma~\ref{lem:qsweight}. 
\end{proof}

	\vspace{3mm}
\noindent{\bf Funding and Competing Interests.}	
We thank the referees for very helpful comments and suggestions which improve the exposition of the paper.

M.L. is partially supported by the National Natural Science Foundation of China (No. 12171333, 1226113149). W.W. and W.Z. are partially supported by the NSF grant DMS-2001351. 

The authors have no competing interests to declare that are relevant to the content of this article.
The authors declare that the data supporting the findings of this study are available within the paper.


\end{document}